\documentclass{gtart_h}

\def\ifplaintex{\expandafter\ifx\csname documentclass\endcsname\relax}

\ifplaintex 
\else
\fi

\def\gt{{\mathsurround=0pt\it $\cal G\mskip-2mu$eometry \&\ 
$\cal T\!\!$opology}}        

\def\gtm{{\mathsurround=0pt\it $\cal G\mskip-2mu$eometry \&\ 
$\cal T\!\!$opology $\cal M\mskip-1mu$onographs}}    

\def\gtp{{\mathsurround=0pt\it $\cal G\mskip-2mu$eometry \&\ 
$\cal T\!\!$opology $\cal P\!$ublications}}  

\def\recd{{\small Received:\qua\receiveddate\ifx\reviseddate\relax
\else\qquad Revised:\qua\reviseddate\fi\par}} 

\def\volumenumber#1{\def\thevolumenumber{#1}}
\def\volumeyear#1{\def\thevolumeyear{#1}}
\def\volumename#1{\def\thevolumename{#1}}
\def\papernumber#1{\def\thepapernumber{#1}}
\def\pagenumbers#1#2{\def\startpage{#1}\def\finishpage{#2}}
\def\published#1{\def\publishdate{#1}}
\def\received#1{\def\receiveddate{#1}}

\def\accepted#1{\def\accepteddate{#1}}

\long\def\asciiabstract#1{\long\def\theasciiabstract{#1}}
\def\asciikeywords#1{\def\theasciikeywords{#1}}

\let\\\par
\let\thevolumenumber\relax\let\thepapernumber\relax
\let\thevolumeyear\relax\let\startpage\relax
\let\finishpage\relax\let\publishdate\relax\let\receiveddate\relax
\let\reviseddate\relax\let\accepteddate\relax\let\theasciititle\relax
\let\theasciiauthors\relax
\let\theasciiabstract\relax\let\theasciikeywords\relax

\let\theerratum\relax\let\theasciiemail\relax
\let\theshortauthors\relax\let\theshorttitle\relax

\def\startpage{1}\def\finishpage{15}\def\thepapernumber{77}

\volumenumber{2}
\volumename{Proceedings of the Kirbyfest}
\volumeyear{1999}

\long\def\maketitlep{   

\count0=\startpage

\gtm\nl        
{\small Volume \thevolumenumber: \thevolumename\nl 
\ifx\theerratum\relax\else Erratum \erratumnumber\nl\fi
Pages \startpage--\finishpage\nl}

\vglue 0.1truein   

{\parskip=0pt\leftskip 0pt plus 1fil\def\\{\par\smallskip}{\ifplaintex\large
\else\Large\fi\bf\thetitle}\par\medskip}   
\vglue 0.05truein 

{\parskip=0pt\leftskip 0pt plus 1fil\def\\{\par}{\sc\theauthors}
\par\medskip}%
 
\vglue 0.03truein 

{\small\leftskip 25pt\rightskip 25pt{\bf Abstract}\stdspace\theabstract

{\bf AMS Classification}\stdspace\theprimaryclass
\ifx\thesecondaryclass\relax\else; \thesecondaryclass\fi\par
{\bf Keywords}\stdspace \thekeywords\par}\vglue 7pt

}   

\font\phead=cmsl9 scaled 950
\font=cmsl9 scaled 1050
\font\pnum=cmbx10 scaled 913
\font\lnum=cmbx10 
\font\pfoot=cmsl9 scaled 950
\font=cmsl9 scaled 1050
\ifplaintex
\headline{\vbox to 0pt{\vskip -4.5mm\line{\small\phead\ifnum
\count0=\startpage ISSN 1464-8997 (on line)
1464-8989 (printed) \hfill {\pnum\folio}\else\ifodd\count0\def\\{ }%
\ifx\theshorttitle\relax\thetitle\else\theshorttitle\fi\hfill{\pnum\folio}
\else\def\\{ and }{\pnum\folio}\hfill\ifx\theshortauthors\relax\theauthors
\else\theshortauthors\fi\fi\fi}\vss}}
\footline{\vbox to 0pt{\vglue 0mm\line{\small\pfoot\ifnum\count0=\startpage
Published \publishdate:\qua\copyright\ \gtp\hfill\else
\gtm, Volume \thevolumenumber\ (\thevolumeyear)\hfill\fi}\vss
}}
\else
\makeatletter
\def\@oddhead{{\small\lhead\ifnum\count0=\startpage ISSN 1464-8997 (on line)
1464-8989 (printed) \hfill {\lnum\number\count0}\else\ifodd\count0
\def\\{ }\ifx\theshorttitle\relax \thetitle \else\theshorttitle\fi\hfill
{\lnum\number\count0}\else\def\\{ and }{\lnum\number\count0}
\hfill\ifx\theshortauthors\relax 
\theauthors\else\theshortauthors\fi\fi\fi}}\def\@evenhead{@oddhead}
\def\@oddfoot{\small\lfoot\ifnum\count0=\startpage Published \publishdate:\qua\copyright\ \gtp\hfill\else
\gtm, Volume \thevolumenumber\ (\thevolumeyear)\hfill\fi}
\def\@evenfoot{@oddfoot}
\makeatother
\fi

\let\maketitlepage\maketitlep
\let\makeshorttitle\maketitlepage
\let\maketitle\maketitlepage

\newwrite\gtoutfile
\long\gdef\makeheadfile{  
{\def\\{, }\def\s{ }
\immediate\openout\gtoutfile head.xxx
\immediate\write\gtoutfile{Proxy-for: \ifx\theasciiauthors\relax
\theauthors\else\theasciiauthors\fi\s<\ifx\theasciiemail\relax\theemail\else\theasciiemail\fi>}
\immediate\write\gtoutfile{\noexpand\\}
\immediate\write\gtoutfile{Authors: \ifx\theasciiauthors\relax
\theauthors\else\theasciiauthors\fi}
{\def\\{ }\immediate\write\gtoutfile{Title: \ifx\theasciititle\relax
\thetitle\else\theasciititle\fi}}
\immediate\write\gtoutfile{Subj-class: GT or SG, GR etc}
\immediate\write\gtoutfile{MSC-class: \theprimaryclass\ifx\thesecondaryclass\relax\else, \thesecondaryclass\fi}
\immediate\write\gtoutfile{Journal-ref: Geom. Topol. Monogr. \thevolumenumber\s
(\thevolumeyear) \startpage-\finishpage}
\immediate\write\gtoutfile{Comments: Published by Geometry and Topology Monographs at}
\immediate\write\gtoutfile{\s\s\s  http://www.maths.warwick.ac.uk/gt/GTMon\thevolumenumber/paper\thepapernumber.abs.html}
\immediate\write\gtoutfile{\noexpand\\}
\immediate\write\gtoutfile{}
\ifx\theasciiabstract\relax
\immediate\write\gtoutfile{\theabstract}\else
\immediate\write\gtoutfile{\theasciiabstract}\fi
\immediate\write\gtoutfile{}
\immediate\write\gtoutfile{\noexpand\\}
\immediate\write\gtoutfile{}
\immediate\closeout\gtoutfile}}  

\def\maketitlepage{\maketitlep\makeheadfile}
\let\makeshorttitle\maketitlepage
\let\maketitle\maketitlepage

\volumenumber{7}
\volumename{Proceedings of the Casson Fest}
\volumeyear{2004}
\papernumber{3}
\pagenumbers{69}{100}
\received{5 August 2003}
\accepted{21 March 2004}
\published{17 September 2004}

\usepackage{amssymb} 
\usepackage{amsmath,amscd} 

\newtheorem{thm}{Theorem}[section]    
\newtheorem{lem}[thm]{Lemma}          
\newtheorem{prop}[thm]{Proposition}
\numberwithin{equation}{section}
\newcommand{\Hcal}{\mathcal{H}}
\newcommand{\C}{\mathbb{C}}
\newcommand{\R}{\mathbb{R}}
\newcommand{\Z}{\mathbb{Z}}
\newcommand{\psl}[2]{\operatorname{PSL}_{#1}(#2)}
\newcommand{\so}[2]{\operatorname{SO}_{#1}(#2)}
\newcommand{\Hom}{\operatorname{Hom}}
\newcommand{\diff}{\operatorname{Diff}}
\newcommand{\Tcal}{\mathcal{T}}
\newcommand{\sym}{\operatorname{Sym}}
\newcommand{\iproduct}[2]{\left\langle#1,#2\right\rangle}
\newcommand{\T}{\mathbb{T}}
\newcommand{\M}{\mathcal{M}}
\newcommand{\met}{\operatorname{Met}}
\newcommand{\Tr}{\operatorname{Tr}}
\newcommand{\tr}{\operatorname{tr}}
\newcommand{\U}{\mathcal{U}}
\newcommand{\re}{\operatorname{re}}
\newcommand{\im}{\operatorname{im}}

\renewcommand{\epsilon}{\varepsilon}

\begin{document}
\title{Minimal surfaces in germs of\\hyperbolic 3--manifolds}
\author{Clifford Henry Taubes}                  
\address{Department of Mathematics, Harvard University\\Cambridge, MA 02138,
USA}
\email{chtaubes@math.harvard.edu}
\begin{abstract}   
This article introduces a universal moduli space for the set whose
archetypal element is a pair that consists of a metric and second
fundamental form from a compact, oriented, positive genus minimal
surface in some hyperbolic $3$--manifold.  This moduli space is a
smooth, finite dimensional manifold with canonical maps to both the
cotangent bundle of the Teichm\"uller space and the space of
$\so{3}{\C}$ representations for the given genus surface.  These two
maps embed the universal moduli space as a Lagrangian submanifold in
the product of the latter two spaces.
\end{abstract}
\asciiabstract{%
This article introduces a universal moduli space for the set whose
archetypal element is a pair that consists of a metric and second
fundamental form from a compact, oriented, positive genus minimal
surface in some hyperbolic 3-manifold.  This moduli space is a smooth,
finite dimensional manifold with canonical maps to both the cotangent
bundle of the Teichmueller space and the space of SO(3,C)
representations for the given genus surface.  These two maps embed the
universal moduli space as a Lagrangian submanifold in the product of
the latter two spaces.}


\asciikeywords{Hyperbolic 3-manifold, minimal surface}                    
\maketitle

\section{Introduction}
Immersed, compact, minimal surfaces are now known to appear in every
compact, hyperbolic $3$--manifold.  Indeed, Schoen and
Yau~\cite{schoen-yau} proved that a Haken hyperbolic manifold has at
least one stable minimal surface.  More recently, Pitts and
Rubinstein~\cite{pr1,pr2,r}, (see also~\cite{colding-delellis})
proved that all compact, hyperbolic 3--manifolds have at least one
unstable, immersed, minimal surface.  Other authors, for example
Freedman, Hass and Scott~\cite{freedman} and Hass and
Scott~\cite{hass} also have foundational papers on this subject.  The
ubiquity of minimal surfaces in hyperbolic $3$--manifolds motivates the
introduction and study of a universal moduli space for the set whose
archetypal element is a pair that consists of a metric and second
fundamental form from a compact, oriented, positive genus minimal
surface in some hyperbolic $3$--manifold.  This article introduces such
a moduli space and takes some (very) small steps to elucidate its
properties.

In this regard, the moduli space introduced below has components that
are labeled in part by the Euler class, $-\chi$, of the surface.  As
explained below, the component with label $\chi$ is a smooth,
orientable manifold of dimension $6\chi$.  Numerologists might notice
that this number is the dimension of the cotangent bundle to the Euler
class $-\chi$ Teichm\"uller space, and that it is also the dimension of
the adjoint action quotient of the space of homomorphisms from the
surface fundamental group to $\psl{2}{\C}$ in its guise as the group
$\so{3}{\C}$.  In fact, the moduli space introduced here admits a
canonical map to each of the latter two spaces, and these maps play a
central role in what follows.

The symbol $\Hcal$ is used below to denote the union of the fixed $\chi$
components of the moduli space.

Before starting, please note that there is a tremendous body of
published research on the subject of minimal surfaces in
$3$--manifolds.  The fact is that such surfaces are an old and well
used tool for studying hyperbolic $3$--manifolds and $3$--manifolds in
general.  Meanwhile, the author pleads the case of a noviate to the
subjects of minimal surfaces and hyperbolic $3$--manifolds.  On the
basis of this meager excuse, the author hereby asks to be forgiven for
his rhyolitic ignorance of the fundamental work of others, and also
for belaboring what may appear obvious to the experts.  

In any event, Rubinstein~\cite{r} has an excellent review of various
aspects of the minimal surface story as applied to $3$--manifolds.
Meanwhile, Colding and Minnicozzi~\cite{cm1} have a recent monograph that
reviews many of the more analytic aspects of the subject.  On the
bigger subject of $3$--manifolds, Scott's exposition ~\cite{scott} is
still very much worth reading.  

Finally, take note that moduli spaces of minimal surfaces have been
introduced by others (Brian White ~\cite{white1} and also Colding and
Minnicozzi~\cite{cm2}, for example); however, the results in these studies
do not appear to speak directly to the moduli space defined here.
	
The remainder of this article is organized as follows: The subsequent
parts of this section provide the precise definition of $\Hcal$ and
describes its canonical maps to the cotangent bundle of Teichm\"uller
space and to 
\begin{equation*}
\Hom(\pi_1(\Sigma); \so{3}{\C})/\so{3}{\C}.
\end{equation*}
Section~\ref{section:structure} explains why $\Hcal$ is a smooth
manifold, while Section~\ref{section:maps} describes the critical loci
of its two canonical maps.  Section~\ref{section:forms} next describes
the pull-backs via these maps of certain natural symplectic structures
on the Teichm\"uller space cotangent bundle and on
$\Hom(\pi_1(\Sigma); \so{3}{\C})/\so{3}{\C}$.
Section~\ref{section:extending} explains why every element in $\Hcal$
arises as a minimal surface in some (usually incomplete) hyperbolic
$3$--manifold.  Finally, Section~\ref{section:neighborhood} describes
an open subset of $\Hcal$ whose elements arise as minimal surfaces in
certain complete, quasi-Fuchsian hyperbolic metrics on
$\Sigma\times\R$.  There are also a number of appendices for novices
that provide derivations of more or less classical formulas.

\subsection{The definition of $\Hcal$}
In what follows, $\Sigma$ denotes a compact, oriented, $2$--dimensional
manifold with negative Euler characteristic.  The absolute value of
the Euler characteristic is denoted as $\chi$.  A pair $(g, m)$ of
Riemannian metric and symmetric section of $T^*\Sigma\otimes
T^*\Sigma$ will be called a ``hyperbolic germ'' on $\Sigma$ if the
following conditions are met:
\begin{subequations}
\begin{equation}\label{eq:1.1a}
\begin{split}
d_Cm_{AB}-d_Bm_{AC} & = 0\\
r + (|m|^2+\frac{1}{3})-k^2 & = 0
\end{split}
\end{equation}
Here, and below, the notation is as follows: First, the subscripts
indicate components with respect to some local frame for $T^*\Sigma$. Second,
$d_C$ denotes the covariant derivative defined by the metric $g$, the norms
are defined by the metric $g$, and the respective functions $r$ and $k$ are
the scalar curvature for the metric $g$ and the trace of $m$ as defined
using the metric $g$.  Finally, repeated indices are to be summed.  

A ``minimal hyperbolic germ'' on $\Sigma$ is a pair, $(g, m)$, of
metric and symmetric tensor that obeys both~\eqref{eq:1.1a} together with
the auxiliary condition
\begin{equation}\label{eq:1.1b}
k = g^{AB}m_{AB}\equiv 0.
\end{equation}
\end{subequations}
With regards to the terminology, an argument is given below to prove
that there is an honest hyperbolic metric on a neighborhood of
$\Sigma\times\{0\}$ in $\Sigma\times\R$ whose respective first and
second fundamental forms on $\Sigma\times\{0\}$ are $g$ and $m$ when
$(g, m)$ is a hyperbolic germ.  The latter metric and that defined by
the line element
\begin{equation}\label{eq:1.2}
ds^2 = (g_{AB} + 2tm_{AB} +
\frac{1}{2}t^2(|m|^2+\frac{1}{3})g_{AB})dz^Adz^B+dt^2 
\end{equation}
agree to order $t^3$ near $t=0$.  In this regard, a metric given
by~\eqref{eq:1.2} has $R_{ij} = -\frac{1}{3}g_{ij}$ at $t=0$ if and
only if~\eqref{eq:1.1a} holds.  Thus, \eqref{eq:1.1a} insures that the
metric in~\eqref{eq:1.2} is hyperbolic to first order at $t = 0$.
Granted~\eqref{eq:1.1a}, the surface $\Sigma\times\{0\}$ has zero mean,
extrinsic curvature with respect to this same metric if and only
if~\eqref{eq:1.1b} holds.  Thus, it is a minimal surface with respect
to the metric in~\eqref{eq:1.2} and to the associated hyperbolic
metric on a neighborhood of $\Sigma\times\{0\}$ in $\Sigma\times\R$.
In this regard, keep in mind that the normalization used here is such
that the scalar curvature of the $3$--dimensional hyperbolic metric is
$-1$; thus its sectional curvatures are $-\frac{1}{3}$.

Define an equivalence relation on the space of hyperbolic germs
whereby any given pair of such germs are identified when one is
obtained from the other by a diffeomorphism of $\Sigma$ that lies in
the component of the identity in $\diff(\Sigma)$.  In this regard, the
infinitesimal form of the diffeomorphism group's action on the space
of hyperbolic germs has a vector field, $v^A\partial_A$, sending
$g_{AB}$ and $m_{AB}$ to
\begin{equation}\label{1.3}
\begin{split}
\delta g_{AB} &= d_Av_B + d_Bv_A\\ 
\text{and}\qquad\delta m_{AB} &= v^Cd_Cm_{AB}+m_{BC}d_Av_C+m_{AC}d_Bv_C.
\end{split}
\end{equation}
Note that this equivalence relation preserves the subset of minimal
hyperbolic germs.  

Let $\Hcal$ denote the quotient.  Thus, $\Hcal = \{\text{minimal
hyperbolic germs}\}/\diff_0(\Sigma)$, where $\diff_0(\Sigma)$ is the
component of the identity in the diffeomorphism group of $\Sigma$.
The set $\Hcal$ inherits the quotient topology with the understanding
that the hyperbolic germs have the induced topology as a subset of the
space of smooth, symmetric, $2\times2$ tensor fields on $\Sigma$.
Please tolerate the notation used below whereby a pair of metric and
traceless, symmetric tensor is said to ``be in'' $\Hcal$.  Of course,
this means that the $\diff_0(\Sigma)$ orbit of the given pair is in
$\Hcal$.

By the way, $\Hcal$ is assuredly non-empty; indeed, if $g$ is a
hyperbolic metric on $\Sigma$ with scalar curvature $-\frac{1}{3}$,
then the pair $(g, 0)$ defines a point in $\Hcal$.  According to the
upcoming Theorem~\ref{thm:2.1}, this is but a small slice of $\Hcal$.
In any event, the space $\Hcal$ is the subject of this article.  What
follows is a brief summary of the story.

The space $\Hcal$ is a smooth, orientable manifold whose dimension is
$6\chi$ where $\chi$ denotes the absolute value of the Euler
characteristic of $\Sigma$.  Moreover, $\Hcal$ admits smooth maps to
the cotangent bundle of $\Sigma$'s Teichm\"uller space and to the
moduli space of flat, $\so{3}{\C}$ connections on $\Sigma$, both
spaces with dimension $6\chi$.  Neither map is proper and both admit
critical points.  This said, here is a surprise: The critical loci of
these maps are identical, this being the loci of pairs in $\Hcal$
where $\Sigma\times\{0\}$ has isotopies in $\Sigma\times\R$ that
preserve its minimality to first order.  However, even as the critical
loci agree, the kernels of the respective differentials are linearly
independent.  The reasons for this coincidence are mysterious, though
almost surely related to the following added surprise: The canonical
symplectic forms on the cotangent bundle to Teichm\"uller space and on
the space of flat $\so{3}{\C}$ connections agree upon pull-back to
$\Hcal$.  In particular, with the signs of these symplectic forms
suitably chosen, these maps immerse $\Hcal$ as a Lagrangian subvariety
in the product of the cotangent bundle to Teichm\"uller space and the
smooth portion of $\Hom(\pi_1; \so{3}{\C})/\so{3}{\C}$.  

By the way, an analysis of the critical loci of these maps from
$\Hcal$ leads to the following observation: The nullity of a compact,
oriented and immersed minimal surface in a hyperbolic $3$--manifold is
no larger than $3$ times the absolute value of its Euler class.

\subsection{$\Hcal$ and Teichm\"uller space}
The additive group of smooth functions acts on the space of metrics to
change the conformal factor.  This is to say that a function $u$ sends
a metric $g$ to $e^{-u}g$.  This action extends to one on the space of
pairs $(g, m)$ with $m$ left unchanged.  Defined in this way, the
first equation in~\eqref{eq:1.1a} is preserved by this action, as is
the equation in~\eqref{eq:1.1b}.  The second equation is not invariant
under such a change.  As demonstrated in the appendix, the second
equation changes as follows:
\begin{equation}
r_{e^{-u}g} + |m|_{e^{-u}g}^2 + \frac{1}{3} =%
 e^u(r_g+\Delta_gu+e^u|m|_g^2 + \frac{1}{3}e^{-u}).
\end{equation}
In any event, the quotient of the space of smooth metrics on $\Sigma$
by the action of the semi-direct product of $\diff_0$ and
$C^{\infty}(\Sigma)$ is called Teichm\"uller space.  Here it is
denoted by $\Tcal$; it is a smooth manifold of dimension $3\chi$.  The
projection from the space of metrics to $\Tcal$ induces a smooth map
from $\Hcal$ to $T\Tcal$.  For certain purposes, it is often
convenient to view this projection as a map to $T^*\Tcal$; this is
done by using the second component of any pair $(g, m)$ to define a
measure on $\Sigma$ with values in $\sym^2(T\Sigma)$, this denoted by
$\hat{m}$.  In particular $\hat{m}$ has components
\begin{equation}
  \hat{m}^{AB} = \det(g)^{\frac{1}{2}}g^{AC}g^{BD}m_{CD}
\end{equation}
in a local coordinate frame.  The tensor-valued measure $\hat{m}$ then
defines a linear functional on the tangent space to the space of
metrics, this the functional whose value on a symmetric tensor $h$ is
given by
\begin{equation}\label{1.6}
  \int_{\Sigma}\hat{m}^{AB}h_{AB} =
  \int_{\Sigma}g^{AC}g^{BD}m_{CD}h_{AB}d\operatorname{vol}_g. 
\end{equation}
Because m is traceless and obeys the top equation in~\eqref{eq:1.1a},
the linear function defined by~\eqref{1.6} annihilates tangent vectors
at $g$ to the orbit of $\diff_0(\Sigma)\times C^{\infty}(\Sigma)$.
Thus, it descends with $g$ to define an element in $T^*\Tcal$.

Granted the preceding, any pair $(g, m)\in\Hcal$ gives a point in
$T^*\Tcal$, and these assignments thus define a canonical map from
$\Hcal$ to $T^*\Tcal$.

\subsection{$\Hcal$ and $\mathrm{SO}_3(\C)$} 
Change gears now to consider the assertion that any given pair $(g,
m)\in\Hcal$ can also be used to define a flat $\so{3}{\C}$ connection
over $\Sigma$.  To elaborate, these will be connections on the
complexification, $E$, of $T^*\Sigma\oplus\R$.  For the purpose of
defining such a connection, use the metric $g$ with the Euclidean
metric to define an inner product, $\langle,\rangle_g$, on
$T^*\Sigma\oplus\R$.  The latter induces a $\C$--bilinear, quadratic form on $E$
which will be denoted in the same way.  Since $E$ is the
complexification of a real $3$--plane bundle, it inherits a tautological
complex conjugation involution; this denoted by an overbar.  Of
course, this complex conjugation and the bilinear form define a
hermitian inner product on $E$ in the usual way; this the polarization
of the norm whose square sends $\eta$ to
$\iproduct{\overline{\eta}}{\eta}_g$.

An $\so{3}{\C}$ connection on $E$ is defined by its covariant
derivative, $\nabla$.  In this regard, the latter must have the property that
\begin{equation}
  d\iproduct{\eta}{\eta'}_g = \iproduct{\nabla\eta}{\eta'}_g+
    \iproduct{\eta}{\nabla\eta'}_g
\end{equation}
whenever $\eta$ and $\eta'$ are sections of $E$.  

The flat $\so{3}{\C}$ connection defined by $(g, m)$ is best expressed
using a local, oriented orthonormal frame $\{e^A\}_{A=1,2}$ for
$T^*\Sigma$ and the unit vector, $e^3$, on $\R$ that points in the positive
direction.  With respect to such a frame, any given section $\eta$ of $E$
that is defined where the local frame is defined can be written as a
column vector,
\begin{equation}
  \eta = \begin{pmatrix} \eta_B\\ \eta_3 \end{pmatrix}
\end{equation}
The covariant derivative defined by the pair $(g, m)$ sends such a
section to $\nabla\eta = \nabla_A\eta e^A$ with
\begin{equation}\label{1.8}
  \nabla_A\eta\equiv \begin{pmatrix} d_A\eta_B + \theta_{AB}\eta_3\\
     d_A\eta_3-\theta_{AC}\eta_C\end{pmatrix}
\end{equation}
and
\begin{equation}\label{1.9}
  \theta_{AB} = m_{AB}+\frac{i}{\sqrt{6}}\varepsilon_{AB}.
\end{equation}
Here, $d_A$ now denotes the metric's covariant derivative when acting
on a section of $T^*\Sigma_{\C}$, and it denotes the exterior
derivative when acting on a complex-valued function.

For use below, note parenthetically that three related, flat
connections can also be defined on $E$.  To define the first,
introduce the isometry $\T\co E\to E$ that acts trivially on the $\R$
factor and as multiplication by $-1$ on the $T^*\Sigma$ factor.  This
done, the new covariant derivative, $\nabla'$, is given by the
formula $\nabla'\equiv\T\nabla\T$.  In particular, $\nabla'$ is given
with respect to the orthonormal frame $\{e^A\}$ by replacing every
$\theta_{AB}$ by $-\theta_{AB}$ in~\eqref{1.8}.  Note that $\T$
defines an automorphism of $E$ with values in $\operatorname{SO}(3)$,
so $\nabla$ and $\nabla'$ are gauge equivalent.  

The remaining two flat connections are also defined by their covariant
derivatives, these denoted by $\overline{\nabla}$ and
$\overline{\nabla}'$.  In this regard,
$\overline{\nabla}'\equiv\T\overline{\nabla}\T$ while
$\overline{\nabla}$ is defined by the formula
\begin{equation}
  \overline{\nabla}\eta\equiv\overline{\nabla\overline{\eta}}.
\end{equation}
In particular, the formula for $\overline{\nabla}$ with respect to a
local frame is obtained from~\eqref{1.8} by replacing every
$\theta_{AB}$ by its complex conjugate.  Note that the formal $L^2$
adjoint of $\nabla$ as defined using the Hermitian inner product on
$C^{\infty}(E)$ is $\nabla^{\dagger}=-\overline{\nabla}$.  Likewise,
$\nabla'^{\dagger} = -\overline{\nabla}'$.

With an inner product fixed on $T*\Sigma\oplus\R$, let $\M$ denote the
moduli space of flat $\so{3}{\C}$ connections on the complexification,
$E$.  This is the quotient of the space of flat $\so{3}{\C}$
connections by the action of the group of automorphisms of $E$.  In
this regard, note that parallel transport around a fixed basis for
$\pi_1(\Sigma)$ identifies $\M$ with the quotient in
$(\so{3}{\C})^{\chi+2}$ of a codimension $6$ subvariety by the adjoint
action of $\so{3}{\C}$.  Here, the subvariety in question is the inverse
image of the identity in $\so{3}{\C}$ under the map that sends a
$(\chi+2)$--tuple of matrices $(U_1,\dots, U_{\chi+2})$ to
\begin{equation}
  (U_1U_2U_1^{-1}U_2^{-1})\dots(U_{\chi+1}U_{\chi+2}U_{\chi+1}^{-1}U_{\chi+2}^{-1}).
\end{equation}
It is useful at times to fix a base point, $z_0\in\Sigma$, and take
the quotient by the group of automorphisms that act as the identity on
the fiber of $E$ over $z_0$.  The choice of a generating basis for the
fundamental group $\pi_1(\Sigma; z_0)$ identifies the latter space,
$\M^0$, with the aforementioned codimension $6$ subvariety in $(\so{3}{\C})^{\chi+2}$. With regards to the structure of $\M^0$, note that the
differential of the map $f$ has a cokernel only at those
$(\chi+2)$--tuples that consist of matrices that all fix a non-zero element
in the lie algebra under the adjoint representation.  A connection
that corresponds to the latter sort of $(\chi+2)$--tuple is said to be
``reducible''.  In particular, note that a connection is reducible if
and only if there exists a non-zero section of $E$ that is annihilated
by the corresponding covariant derivative.  The complement in $\M$ of the
set reducible connections is a smooth manifold of dimension $6\chi$.

The association of the covariant derivative $\nabla$ to a pair $(g, m)\in
\Hcal$ defines a continuous map from $\Hcal$ to $\M$.  As argued below,
the image of $\Hcal$ avoids the reducible connections.

\section{The structure of $\Hcal$}\label{section:structure}
The following theorem is the principle result of this subsection:

\begin{thm}\label{thm:2.1}
 The space $\Hcal$ has the structure of a smooth, orientable manifold
 of dimension $6\chi$.  Moreover, $\Hcal$ comes equipped with a smooth
 action of $S^1$, this provided by the map that sends $\tau\in S^1 =
 \R/(2\pi\Z)$ and $(g,m)\in\Hcal$ to $(g, \cos\tau\,m + \sin\tau\,
 \varepsilon\cdot m)$ where $\varepsilon\cdot m$ is the symmetric,
 traceless tensor with components $\varepsilon_{AC}m_{CB}$ in a local,
 $g$--orthonormal frame.
\end{thm}

Note that this circle action is free away from the locus in
$\Hcal$ whose elements are pairs of the form $(g, 0)$ where $g$ is a metric
on $\Sigma$ with constant scalar curvature $-\frac{1}{3}$.  

The remainder of this subsection is occupied with the proof of this
proposition.

To start the proof, remark first that the linearization of the
equations in~\eqref{eq:1.1a} and~\eqref{eq:1.1b} about any given pair
$(g, m)$ defines a differential operator; the latter denoted by
$L_{(g,m)}$ in what follows.  In local coordinates, this operator
sends a pair consisting of a symmetric tensor and a symmetric,
traceless tensor to a pair consisting of a $1$--form and a functions.
In particular, the map sends a pair $(h_{AB},n_{AB})$ to the
vector/function pair
\begin{equation}\label{2.1}
  \begin{split}
    \gamma_B \equiv\; &
    \epsilon_{CA}d_Cn_{AB}+\frac{1}{2}\epsilon_{AC}(d_Ch_{DB})m_{AD}\\
     & +\frac{1}{2}(\epsilon_{EF}d_Eh_{FC})m_{BC}+\frac{1}{2}\epsilon_{CB}d_C(h_{EF}m_{EF}),
    \\
  \gamma_3 \equiv\; &
  \frac{1}{2}(\frac{1}{3}-|m|^2)h_{AA}+d_Ad_Bh_{AB}-d_Bd_Bh_{AA}+2n_{AB}m_{AB}. 
  \end{split}
\end{equation}
In this regard, the first order variation of $m_{AB}$ is not $n_{AB}$, but
rather $n_{AB} + \frac{1}{2}g_{AB}(h_{EF}m_{EF})$; this an imposition from~\eqref{eq:1.1b}.

This $L_{(g,m)}$ extends as a bounded, linear map from various Sobolev
space completions of its domain to corresponding completions of the
range space as a semi-Fredholm map, a map with infinite dimensional
kernel, closed range and finite dimensional cokernel.  For example, it
has such an extension from $L^2_2(\sym^2 T^*\Sigma)\oplus L^2_1(\sym^2
T^*\Sigma)$ to $L^2(T^*\Sigma)\oplus L^2(\Sigma)$.

By the way, the kernel is infinite dimensional due to the fact that
all pairs $(h, n)$ that induce~\eqref{1.3} are in its kernel.  This said,
introduce the operator $l_{(g,m)}$ that maps $1$--forms to the domain of
$L_{(g,m)}$ by the rule
\begin{equation}\label{2.2}
  v_B\rightarrow (d_Av_B+d_Bv_A,
  v^Cd_Cm_{AB}+m_{BC}d_Av_C+m_{AC}d_Bv_C-g_{AB}(m_{EF}d_Ev_F)).
\end{equation}
The restriction of $L_{(g,m)}$ to the $L^2$--orthogonal complement of the image
of $l_{(g,m)}$ is then Fredholm.  Let $L_{*(g,m)}$ denote this restricted
operator.  The implicit function theorem in conjunction with ``off the
shelf'' differential equation technology can be employed to prove the
following:

\begin{lem}
Let $(g, m)\in\Hcal$.  There exists a ball $B$ about the origin in
the kernel of $L_{*(g,m)}$, a smooth map,
\begin{equation}
  f\co B\to\operatorname{cokernel}(L_{*(g,m)})
\end{equation}
that maps $0$ to $0$, and a homeomorphism from $f^{-1}(0)$ to a
neighborhood of $(g, m)$ in $\Hcal$ that sends the origin to $(g,
m)$. Moreover, if $L_{*(g,m)}$ has trivial cokernel at a given $(g,
m)$, then the corresponding cokernel is trivial at all points in some
neighborhood of $(g, m)$ in $\Hcal$; and this neighborhood has the
structure of a smooth manifold of dimension $6\chi$.
\end{lem}

Those points in $\Hcal$ where the cokernel of the operator
$L_{(\cdot)}$ is trivial will be called ``regular points''.  This
understood, the fact that the space $\Hcal$ is a manifold follows from

\begin{prop}\label{prop:2.3}
  All points in $\Hcal$ are regular points.
\end{prop}

\begin{proof}[Proof of Proposition~\ref{prop:2.3}]
\renewcommand{\qedsymbol}{}
A pair $(\sigma_B,\sigma_3)$ is in the cokernel of $L_{(g,m)}$ if and
only if it is $L^2$--orthogonal to all pairs $(\gamma_B,\gamma_3)$ that
can be written as in~\eqref{2.1}. It proves useful for this and other
purposes to replace $(\sigma_B,\sigma_3)$ with
$v_b\equiv-\epsilon_{BC}\sigma_C$ and $v_3\equiv-2\sigma_3$.  As is
explained in Appendix~\ref{appendix:D}, the fact that
$(\sigma_B,\sigma_3)$ is in the cokernel of $L_{(g,m)}$ implies that
\begin{equation}\label{2.5}
  \begin{split}
    \eta_B&\equiv
    v_B+i\sqrt{6}(-\epsilon_{CB}d_Cv_3+v_E\epsilon_{CB}m_{EC}),\\
    \eta_3&\equiv v_3+i\sqrt{\frac{3}{2}}\epsilon_{EF}d_Ev_F,\\
    u_B&\equiv 0,\quad\text{and}\\
    u_3&= -i\sqrt{\frac{3}{2}}d_Cv_C
  \end{split}
\end{equation}
gives respective complex-valued sections $\eta$ and $u$ of
$T^*\Sigma\oplus\R$ that obey the equation
\begin{equation}\label{2.6}
  \nabla_A\eta = \epsilon_{AB}\nabla_B'u.
\end{equation}
Granted that such is the case, it then follows that both
\begin{equation}
  \nabla_A'\nabla_A\eta = 0\quad\text{and}\quad\nabla_A\nabla_A'u = 0.
\end{equation}
With (2.7) noted, the proof of Proposition 2.3 is obtained by invoking
\end{proof}

\begin{lem}\label{lem:2.4}
 The operators $\nabla_A'\nabla_A$ and $\nabla_A\nabla_A'$ are equal.
 Moreover, if $\eta$ is annihilated by either, then $\eta_B=0$ and
 $\eta_3$ obeys the equation
 \begin{equation}
   d_Ad_A\eta_3+(|m|^2-\frac{1}{3})\eta_3 = 0.
 \end{equation}
\end{lem}

The proof of the Lemma~\ref{lem:2.4} is given momentarily.  To obtain
Proposition~\ref{prop:2.3}, note that in the case at hand, with $\eta$
given by~\eqref{2.5}, this lemma implies that the cokernel element to
$L_{(g,m)}$ defined by $(v_B, v_3)$ has $v_B\equiv 0$ since the latter
is the real part of $\eta_B$.  This then implies via the third and
fourth lines in~\eqref{2.5} that $u\equiv 0$ and then~\eqref{2.6}
requires that $v_3 = 0$ as well.

\begin{proof}[Proof of Lemma~\ref{lem:2.4}]
Written in their components, the equation $\nabla_A'\nabla_A\eta
= 0$ and the equation $\nabla_A\nabla_A'u=0$ read:
\begin{equation}
  \begin{split}
    d_Ad_A\eta_B+\frac{1}{2}(|m|^2-\frac{1}{3})\eta_B+\frac{i}{\sqrt{6}}(m_{AB}\epsilon_{AC}+\epsilon_{AB}m_{AC})\eta_C
    & = 0\\
    d_Ad_A\eta_3+(|m|^2-\frac{1}{3})\eta_3&= 0.
  \end{split}
\end{equation}
To see that only $\eta_B = 0$ solves the top equation, take its
$\C$--linear inner product with $\overline{\eta_B}$ and then integrate
over $\Sigma$.  After an integration by parts, one obtains the
equality
\begin{equation}\label{2.10}
  \int_{\Sigma}(d_A\overline{\eta_B}d_A\eta_B) =%
  \frac{1}{2}\int_{\Sigma}(|m|^2-\frac{1}{3})\overline{\eta_B}\eta_B+%
    \frac{i}{\sqrt{6}}\int_{\Sigma}(m_{AB}\epsilon_{AC}+\epsilon_{AB}m_{AC})\overline{\eta_B}\eta_C. 
\end{equation}
Now, the second term on the right side in~\eqref{2.10} is purely
imaginary by virtue of the fact that the tensor $m_{AB}\epsilon_{AC}$
is symmetric when~\eqref{eq:1.1b} holds.  As the other terms in~\eqref{2.10}
are real, it follows that
\begin{equation}
  \int_{\Sigma}(d_A\overline{\eta_B}d_A\eta_B)=%
   \frac{1}{2}\int_{\Sigma}(|m|^2-\frac{1}{3})\overline{\eta_B}\eta_B.
\end{equation}
To see that no such equality can hold if $\eta_B$ is nonzero, note that
\begin{equation}
  \int_{\Sigma}(d_A\overline{\eta_B}d_A\eta_B) = \int_{\Sigma}(|\epsilon_{AC}d_A\eta_C|^2+|d_A\eta_A|^2)-\frac{1}{2}\int_{\Sigma}r\overline{\eta_B}\eta_B.
\end{equation}
In particular, with the value for $r$ from~\eqref{eq:1.1a}, these last
two equations imply that
\begin{equation}
  \frac{1}{2}\int_{\Sigma}(|m|^2+\frac{1}{3})\overline{\eta_B}\eta_B\leq\frac{1}{2}\int_{\Sigma}(|m|^2-\frac{1}{3})\overline{\eta_B}\eta_B
\end{equation}
which requires that $\eta_B\equiv 0$.
\end{proof}

As is explained next, the fact that $\Hcal$ is orientable follows
using the Atiyah--Singer index theorem for families of operators.  To
elaborate, remark first that the family in question is that defined by
the pair of operators $D_{(g,m)}\equiv (L_{(g,m)}, l_{(g,m)}^*)$,
where $l_{(g,m)}^*$ is the formal $L^2$--adjoint of the operator
$l_{(g,m)}$ that appears in~\eqref{2.2}.  Thus, $D_{(g,m)}$ maps a
pair consisting of a symmetric tensor and a symmetric, traceless
tensor to a triple consisting of a $1$--form, a function, and a vector
field on $\Sigma$.  This operator is defined for any pair $(g, m)$ of
metric on $\Sigma$ and traceless, symmetric tensor.  It is equivariant
with respect to the action of $\diff_0(\Sigma)$ on the space of such
pairs, and it is Fredholm for any given pair $(g, m)$.  This noted,
the Atiyah--Singer index theorem for families of
operators~\cite{atiyah-singer} asserts that the kernels and cokernels
of this operator as the pair $(g, m)$ vary defines a class in the
real, $\diff_0(\Sigma)$--invariant real $K$--theory of the space of
pairs of metric and symmetric, traceless tensor on $\Sigma$.

The relevance of this last point to the orientation question stems
from the fact that the kernel of $D_{(g,m)}$ along the inverse image
of $\Hcal$ projects to $\Hcal$ as its tangent space.  Thus, the first
Stiefel--Whitney class of the corresponding $K$--theory class along
$\Hcal$ is that of the tangent bundle to $\Hcal$.  In particular,
$\Hcal$ is orientable if the first Stiefel--Whitney class of this
$K$--theory class is zero.

To prove that such is the case, note that the space of pairs
$(g, m)$ deformation retracts in the obvious manner onto its $m = 0$
subspace.  On the latter subspace, a given version of this operator
sends a pair $(h, n)$ in its domain to
\begin{multline}
  (\epsilon_{AC}d_Cn_{AB},
  \frac{1}{6}h_{AA}+d_Ad_Bh_{AB}-d_Bd_Bh_{AA},-2d_Ah_{AB})\\ \in
  C^{\infty}(T^*\Sigma)\oplus C^{\infty}(\Sigma)\oplus C^{\infty}(T\Sigma).
\end{multline}
The kernel of such a $(g, 0)$ version of $D_{(g,m)}$ consists of the pairs
$(h, n)$ where both are traceless and symmetric and both the $w=h$ and $w = n$ versions of the equation $\epsilon_{AC}d_Cw_{CB} = 0$ are satisfied.
Meanwhile, the cokernel is trivial in all cases because no metric on a
surface with genus greater then $1$ has a non-trivial, conformal Killing
vector field.  

As just identified, the kernel space of $D_{(g,0)}$ admits the evident
almost complex structure that sends a pair $(h, n)$ to $(\epsilon\cdot
h, \epsilon\cdot n)$; thus it can be viewed as a vector space over
$\C$.  This complex structure is compatible with the $\diff_0(\Sigma)$
action, and so the $\diff_0(\Sigma)$--equivariant $K$--theory class in
question comes from a complex $K$--theory class.  As such, its first
Stiefel--Whitney class is automatically zero.

\section{The maps from $\Hcal$ to $\M$ and $T^*\Tcal$}\label{section:maps}
The purpose of this section is to describe certain aspects of the maps
from $\Hcal$ to $\M$ and from $\Hcal$ to $T^*\Tcal$.  The map to $\M$ is
considered first.  In this regard, note that Lemma~\ref{lem:2.4} has the
following immediate corollary:

\begin{prop}\label{prop:3.1}
 The image of $\Hcal$ in $\M$ has no reducible connections. 
\end{prop}

\begin{proof}[Proof of Proposition~\ref{prop:3.1}]
Indeed, if $A$ is reducible, then there is a non-trivial section of $E$
that is annihilated by $\nabla_A$.  Such a section is then annihilated by
$\nabla_A'\nabla_A$ and so it has only its third component non-trivial.
However, since $\theta_{AB}$ is invertible, this is impossible if it is
annihilated by $\nabla_A$.
\end{proof}

To study the differentiable aspects of the map to $\M$, a preliminary
digression is needed to introduce a ``Zariski tangent space'' to the
$\so{3}{\C}$--orbit of a flat connection: If $\nabla$ is the corresponding
covariant derivative on $E$, then this Zariski tangent space is the
vector space of $E$--valued $1$--forms that obey
\begin{equation}
  \epsilon_{AB}\nabla_Aw_B=0\quad\text{and}\quad -\overline{\nabla}_Aw_A=0.
\end{equation}
Use $T\M|_{\nabla}$ to denote the vector space of such forms. In the
case that the connection is irreducible, then $T\M|_{\nabla}$ has
dimension $6\chi$.  Note that when a is any given $E$--valued $1$--form,
then Hodge theory provides its decomposition as
\begin{equation}
  a_A = \nabla_A\alpha+\epsilon_{AB}\overline{\nabla}_B\beta+w_A
\end{equation}
where $w\in T\M|_{\nabla}$ and both $\alpha$ and $\beta$ are sections
of $E$.  Note that this decomposition is orthogonal with respect to
the $L^2$--hermitian inner product on $T^*\Sigma\otimes E$.
Furthermore, $w$ is uniquely defined by $v$; and so are $\alpha$ and
$\beta$ when the connection is irreducible.

Now, the tangent space to a pair $(g, m)\in\Hcal$ can be identified
with the vector space of pairs $(h, n)$ that make $(\gamma_B,
\gamma_3)$ from~\eqref{2.1} vanish identically and are such that
\begin{equation}\label{3.3}
  -d_Ah_{AC}+\frac{1}{2}n_{AB}d_Cm_{AB}-d_A(n_{AB}m_{BC}+\frac{1}{2}h_{EF}m_{EF}m_{AC})
  = 0
\end{equation}
In this regard, remember that both $h$ and $n$ are symmetric
tensors with $n$ traceless, and that they define respective first
order deformations of $g$ and $m$ that are given by $\delta g_{AB} =
h_{AB}$ and $\delta m_{AB} = n_{AB} + h_{EF}m_{EF}g_{AB}$.  The
condition in~\eqref{3.3} asserts only that this version of $\delta g$
and $\delta m$ is $L^2$--orthogonal to all versions of $\delta g$ and
$\delta m$ that are given by~\eqref{1.3}.

To proceed, define an $E$--valued $1$--form $v$ by writing its components
$(v_{AB}, v_{A3})$ as
\begin{equation}\label{3.4}
  \begin{split}
    v_{AB} &=
    -\epsilon_{AE}n_{EB}+\frac{1}{2}\epsilon_{AB}h_{CD}m_{CD}-\frac{1}{2}\epsilon_{CB}m_{FA}h_{CF}-\frac{i}{2\sqrt{6}}h_{AB}\\
   v_{A3}&=\epsilon_{AE}(\frac{1}{2}d_Ch_{CE}-\frac{1}{2}d_Eh_{CC}).
  \end{split}
\end{equation}
The vanishing of the right hand side of the equations in~\eqref{2.1} is then
the assertion that
\begin{equation}
  \epsilon_{AC}\nabla_Av_C=0.
\end{equation}
This section $v$ of $E\otimes T^*\Sigma$ is introduced for the
following reason: The ``Zariski'' tangent space to $\M$ at the
connection that defines $\nabla$ is the vector space of sections of
$E\otimes T^*\Sigma$ that obey the equations $\nabla_A^{\dagger}q_A =
0 = \epsilon_{AB}\nabla_Aq_B$.  Moreover, the $L^2$--orthogonal projection
of $v$ onto this space provides the image of the tangent vector $(h, n)$
to $\Hcal$ in this Zariski tangent space.  This understood, the kernel
of the differential at $(g, m)$ of the map from $\Hcal$ to $\M$ is
isomorphic to the subspace of pairs $(h, n)$ for which the corresponding
$v$ in~\eqref{3.4} can be written as $v_A=\nabla_A u$ with $u$ a smooth section of $E$.

\begin{prop}\label{prop:3.2}
The kernel of the differential at $(g, m)$ of the map from $\Hcal$ to
$\M$ is canonically isomorphic to the kernel of the operator
$\Delta_{(g,m)}\equiv -d_Cd_C -|m|^2+\frac{1}{3}$.  Moreover, the
isomorphism in question sends a function
$\sigma\in\operatorname{kernel}(\Delta_{(g,m)})$ to the pair
$(h,n)\in T\Hcal_{(g,m)}$ that are given in a local frame by
\begin{equation}\label{3.6}
  \begin{split}
    h_{AB} =& -(d_A\sigma_B+d_B\sigma_A)-2m_{AB}\sigma \\
    n_{CB}
    =&-(\sigma_Fd_Fm_{BC}+m_{CF}d_B\sigma_F+m_{BF}d_C\sigma_F-g_{CB}m_{AF}d_A\sigma_F)\\
    &+(d_Bd_C\sigma-\frac{1}{2}g_{BC}d_Ad_A\sigma);
  \end{split}
\end{equation}
here $\{\sigma_A\}$ are determined by requiring $(h, n)$ to be
orthogonal to the tangent space to the orbit of $\diff_0(\Sigma)$ through $(g,
m)$.
\end{prop}

\begin{proof}[Proof of Proposition~\ref{prop:3.2}]
There are two parts to the proof.  The first argues that any pair $(h,
n)$ as given by the proposition defines a non-trivial tangent vector to
$\Hcal$ at $(g, m)$.  The second part proves that the kernel of the
differential at $(g, m)$ of the map to $\M$ has the asserted form.

\textbf{Part 1}: Suppose that $d_Cd_C\sigma + (|m|^2 -
\frac{1}{3})\sigma = 0$ and $\sigma$ is not identically zero.  To
prove that the pair $(h, n)$ in~\eqref{3.6} defines a non-zero tangent
vector to $\Hcal$ it is necessary to explain why there is no choice
for $\{\sigma_A\}$ that makes $(h, n)$ identically zero, and why the
pair $(h, n)$ makes $(\gamma_B,\gamma_3)$ in~\eqref{2.1} equal to
zero.  For the first task, note that if $h_{AB}= 0$, then $d_C\sigma_C
= 0$ and also $-|m|^2 u = m_{AB}d_A\sigma_B$.  Thus, the vanishing
of~\eqref{3.6} is equivalent to the equations
\begin{equation}
  \begin{split}
    d_C\sigma_C & = 0,\\
    d_A\sigma_B+d_B\sigma_A-g_{AB}d_C\sigma_C-2m_{AB}\sigma & = 0,\\
    \frac{1}{2}(\frac{1}{3}+|m|^2)g_{AB}-d_Ad_B\sigma+\sigma_Fd_Fm_{AB}+m_{CF}d_B\sigma_F+m_{BF}d_C\sigma_F
    & = 0.
  \end{split}
\end{equation}
With the substitute notation $\sigma_A\rightarrow v_A$ and
$\sigma\rightarrow v_3$, this is a version of~\eqref{D.2} in
Appendix~\ref{appendix:D}.  As established in
Appendix~\ref{appendix:D}, the latter equation is equivalent to the
assertion that complex-valued section, $\eta$, of $T^*\Sigma\oplus\R$
that is given by the top line of~\eqref{2.5} obeys the $u\equiv0$
version of~\eqref{2.6}.  However, according to
Proposition~\ref{prop:3.1}, this means that $\eta\equiv0$, and such can not
be the case unless $\sigma\equiv0$ as well.

The proof that $(h, n)$ makes $(\gamma_A, \gamma_3)$ vanish is left as
a chore for the reader.

\textbf{Part 2}: If $(h, n)$ is mapped to zero by the differential of
the map to $\M$, then $v$ given by~\eqref{3.4} as the form $v = \nabla
u$.  This assumed, write $u\equiv\alpha+i\beta$ where $\alpha$ and
$\beta$ have purely real components.  The equation $v = \nabla u$
implies that
\begin{equation}\label{3.8}
  h_{AB}=-\sqrt{6}(d_A\beta_B+d_B\beta_A+\epsilon_{AB}(\epsilon_{EF}d_E\beta_F+\sqrt{\frac{2}{3}}\alpha_3)
  + 2m_{AB}\beta_3). 
\end{equation}
Note that this last equation implies straight off that 
\begin{equation}
  \epsilon_{EF}d_E\beta_F+\frac{\sqrt{2}}{\sqrt{3}}\alpha_3 = 0
\end{equation}
because $h_{AB}$ is symmetric. Thus, the first line in~\eqref{3.6} is seen
to hold using $\sigma_B = -\sqrt{6}\beta_B$ and $\sigma\equiv
-\sqrt{6}\beta_3$. 

To obtain the required equation of $\sigma = -\sqrt{6}\beta_3$, note
first that as $v_{A3}$ is purely real, the equation $v = \nabla u$
also implies that
\begin{equation}\label{3.10}
  d_A\beta_3-m_{AB}\beta_B-\frac{1}{\sqrt{6}}\epsilon_{AB}\alpha_B = 0.
\end{equation}
To continue, act by $d_A$ on both sides of this last equation, sum
over the index $A$, and invoke the requirement from~\eqref{3.8} that
$m_{AB}h_{AB} = -2\sqrt{6}m_{AB}d_A\beta_B - 2\sqrt{6}|m|^2\beta_3$.
The resulting equation reads
\begin{equation}\label{3.11}
  d_Ad_A\beta_3+|m|^2\beta_3+\frac{1}{2\sqrt{6}}m_{AB}h_{AB}-\frac{1}{\sqrt{6}}\epsilon_{AB}d_A\alpha_B
  = 0.
\end{equation}
Meanwhile, a return to~\eqref{3.4} for $v_{AB}$ finds that 
\begin{equation}\label{3.12}
  m_{AB}h_{AB} = 2\epsilon_{AB}v_{AB} = 2\epsilon_{AB}d_A\alpha_B-\frac{4}{\sqrt{6}}\beta_3
\end{equation}
and so~\eqref{3.11} asserts that 
\begin{equation}
  d_Ad_A\beta_3+(|m|^2-\frac{1}{3})\beta_3 = 0.
\end{equation}
To obtain the second line of~\eqref{3.6}, note that under the given
assumptions, \eqref{3.4} implies that
\begin{equation}\label{3.14}
  n_{CB}-\frac{1}{2}m_{FB}h_{CF} = \epsilon_{CA}d_A\alpha_B+\epsilon_{CA}m_{AB}\alpha_3+\frac{1}{\sqrt{6}}g_{CB}\beta_3.
\end{equation}
Using~\eqref{3.10} to write
$\alpha_B=\sqrt{6}(\epsilon_{DB}d_D\beta_3-\epsilon_{DB}m_{DC}\beta_C)$
then finds
\begin{multline}
  n_{CB}=\frac{1}{2}m_{FB}h_{CF}+\sqrt{6}\epsilon_{CA}d_A(\epsilon_{DB}d_D\beta_3-\epsilon_{DB}m_{DF}\beta_F)\\
  -\sqrt{\frac{3}{2}}\epsilon_{CA}m_{AB}\epsilon_{EF}d_E\beta_F+\frac{1}{\sqrt{6}}g_{CB}\beta_3.
\end{multline}
This then rearranges as
\begin{multline}
 n_{CB}=\frac{1}{2}m_{FB}h_{CF}+\sqrt{6}d_Bd_C\beta_3-\sqrt{6}g_{CB}d_Ad_A\beta_3\\
+\frac{1}{\sqrt{6}}g_{CB}\beta_3-\sqrt{6}d_B(m_{CF}\beta_F)+\sqrt{6}g_{CB}m_{AF}d_A\beta_F\\
-\sqrt{\frac{3}{2}}m_{AB}d_C\beta_A+\sqrt{\frac{3}{2}}m_{AB}d_A\beta_C.
\end{multline}
Meanwhile, \eqref{3.8} also has
\begin{equation}
-\frac{1}{2}m_{FB}h_{CF} =
\frac{1}{2}\sqrt{6}(m_{AB}d_A\beta_C+m_{AB}d_C\beta_A) + \frac{1}{2}\sqrt{6}g_{CB}|m|^2\beta_3.  
\end{equation}
Inserting this into~\eqref{3.12} finds the second assertion in~\eqref{3.6}.
\end{proof}

The next proposition provides a characterization of the critical set
for the map from $\Hcal$ to $T^*\Tcal$.

\begin{prop}\label{prop:3.3}
The kernel of the differential at $(g, m)$ of the map from $\Hcal$ to
$T^*\Tcal$ is canonically isomorphic to the kernel of the operator
$\Delta_{(g,m)}\equiv-d_Cd_C-|m|^2+\frac{1}{3}$.  Moreover, the
isomorphism in question sends a function $\sigma\in\operatorname{kernel}(\Delta_{(g,m)})$ to the
pair $(h, n)\in T\Hcal_{(g,m)}$ that are given in a local frame by
\begin{equation}
  \begin{split}
    h_{AB} & = -(d_A\sigma_B+d_B\sigma_A) + g_{AB}\sigma\\
    n_{CB} & = -(\sigma_Fd_Fm_{BC}+m_{CF}d_B\sigma_F+m_{BF}d_C\sigma_F-g_{CB}m_{AF}d_A\sigma_F)
  \end{split}
\end{equation}
where $\{\sigma_A\}$ are determined by requiring $(h, n$) to be
orthogonal to the tangent space to the orbit of $\diff_0(\Sigma)$
through $(g, m)$.
\end{prop}

\begin{proof}[Proof of Proposition~\ref{prop:3.3}]
Any element in the kernel of the differential of the map to $T^*\Tcal$
must have the form given above for some function $\sigma$, so the
proposition follows by verifying two assertions: First, if $\sigma$ is
a non-trivial element of the kernel of $\Delta_{(g,m)}$, then no
choice for $\{\sigma_A\}$ makes~\eqref{3.6} vanish. Second, if $(h,
n)\in T\Hcal|_{(g,m)}$ and in the kernel of the differential of the
map to $T*\Tcal$, then $\sigma\in\operatorname{kernel}(\Delta_{(g,m)})$.

To verify the first assertion, suppose that the right hand side
of~\eqref{3.14} is zero.  Now define a real-valued section $(v_A,
v_3)$ of $T^*\Sigma\oplus\R$ by setting $v_A\equiv\sigma_A$ and
$v_3\equiv 0$.  This section then satisfies~\eqref{D.6} with $\kappa =
0$ and $\sigma = d_Cv_C$.  Now, as a consequence of
Lemma~\ref{lem:D.1}, the section $\eta$  as defined in~\eqref{D.7} obeys
$\nabla_A'\nabla_A\eta = 0$ and then Lemma~\ref{lem:2.4} asserts that $\eta_B =
0$.  Thus, $v_B = \sigma_B = 0$ and so $\sigma = 0$.

To verify the second assertion, it is sufficient to verify that the
pair $(h, n)$ in~\eqref{3.14} make~\eqref{2.1}'s pair $(\gamma_B,
\gamma_3)$ vanish if and only if
$\sigma\in\operatorname{kernel}(\Delta_{(g,m)})$.  In this regard, it
is sufficient to consider the case where $\{\beta_A\equiv 0\}$.  In
this case, $n_{CB}\equiv 0$ and $\gamma_B$ vanishes identically and
the second line of~\eqref{2.1} asserts that $\sigma\in\operatorname{kernel}(\Delta_{(g,m)})$.
\end{proof}

It is almost surely the case that neither the map from $\Hcal$ to
$T^*\Tcal$ nor that from $\Hcal$ to $\M$ is proper.  In this regard,
known results about the behavior of sequences of immersed surfaces in
a fixed 3--manifold (see, e.g.~\cite{anderson,white2,schoen2}) are surely
relevant.  In any event, non-convergent sequences in $\Hcal$ can be
characterized to some extent.

\section{$\Hcal$ and symplectic forms on $T^*\Tcal$ and $\M$}%
\label{section:forms}
Propositions~\ref{prop:3.2} and~\ref{prop:3.3} nix certain obvious
candidate for a symplectic form on $\Hcal$.  To elaborate, let
$\met(\Sigma)$ denote the space of smooth Riemannian metrics on $\Sigma$ and
let $T^*\met(\Sigma)$ denote the bundle over $\met(\Sigma)$ whose
fiber is the space of symmetric, measured valued sections of
$\Lambda^2T\Sigma$.  Note that this space is a linear subspace of the
honest cotangent bundle of $\Sigma$, the latter being a space of
distribution valued sections of $\Lambda^2T^*\Sigma$.  In any event,
$T^*\met(\Sigma)$ has a canonical, $\diff_0$--invariant symplectic
form; it is defined as follows: If $(\delta g, \delta\hat{m})$ and
$(\delta g', \delta\hat{m}')$ are tangent vectors to $T^*\met$, then their
symplectic pairing is by definition
\begin{equation}\label{4.1}
  \int_{\Sigma}(\delta g_{AB}\delta\hat{m}^{'AB} - \delta
  g_{AB}'\delta\hat{m}^{AB}). 
\end{equation}
Meanwhile, the set of pairs $(g, m)$ that obey the conditions in (1.1)
embeds in $T^*\Sigma$ by the map that sends any given $(g, m)$ to the
pair $(g, \hat{m})$ where $\hat{m}$ is given in local coordinates by
the tensor
$\hat{m}^{AB}\equiv\det(g)^{\frac{1}{2}}g^{AC}g^{BD}m_{CD}$.  Thus,
$\Hcal$ embeds in the quotient of $T^*\met(\Sigma)$ by the action of
$\diff_0(\Sigma)$.  In this regard, the symplectic form in~\eqref{4.1}
does not descend to the whole of the latter.  However, it does
formally descend to the $\diff_0(\Sigma)$ quotient of the zero set in
$T^*\met(\Sigma)$ of a certain ``moment map''.  This map, denoted by
$\wp$, is the map from $T^*\met(\Sigma)$ to the dual of the space of
smooth $1$--forms on $\Sigma$ that sends a pair $(g, \hat{m})$ to the
linear functional on $C^{\infty}(T\Sigma)$ that assigns to any given
vector field $(v_A)$ the value
\begin{equation}\label{4.2}
  \wp(g, \hat{m})\cdot v\equiv\int_{\Sigma}\hat{m}^{AB}d_Av_B.
\end{equation}
As the zero set of $\wp$ consists of those pairs $(g,\hat{m})$ such
that $d_A\hat{m}^{AB} = 0$, it follows from the first line
in~\eqref{eq:1.1a} and~\eqref{eq:1.1b} that $\Hcal$ lies in
$\wp^{-1}(0)/\diff_0(\Sigma)$.  As a consequence, the pairing in~\eqref{4.1}
defines a closed $2$--form on $\Hcal$.  This form is denoted here by $\omega_{\Hcal}$.

To be explicit, the form $\omega_{\Hcal}$ is defined on any two pair
of tangent vectors at the same given point in $\Hcal$ as follows: Take
any pair $(g, m)$ that projects to the given point in $\Hcal$ and
choose respective pairs $(h, n)$ and $(h', n')$ that project to
the given tangent vectors.  Then, the pairing of the two vectors down
on $\Hcal$ with the closed form is
\begin{equation}\label{4.3}
  \int_{\Sigma}(h_{AB}n_{AB}'-h_{AB}'n_{AB}).
\end{equation}
As~\eqref{4.1} defines a closed form on $T^*\met(\Sigma)$,
so~\eqref{4.3} defines a closed $2$--form on $\Hcal$.  Thus, the latter
is symplectic on $\Hcal$ if it is non-degenerate.  However, as will
now be explained, it is degenerate precisely on the critical locus for
the map to $T^*\Tcal$.  To see why, note first that the additive group
of smooth functions acts on $T^*\met(\Sigma)$ where by a function
$\sigma$ sends a given pair $(g, \hat{m})$ to $(e^{-\sigma}g,
e^{\sigma}\hat{m})$.  This action is compatible with that of $\diff_0$
if the joint action is viewed as one of a semi-direct product group.
In any event, the action of $C^{\infty}(\Sigma)$ also preserves the
symplectic form, it preserves $\wp^{-1}(0)$, and it has a moment map,
$p$, whose zero set is $\diff_0(\Sigma)$ invariant.  In particular,
the latter consists of the pairs $(g, \hat{m})$ for which
$g_{AB}\hat{m}^{AB} = 0$.  In this regard, note that any pair $(g, m)$
that solves (1.1) maps to the subset $\wp^{-1}(0)\cap p^{-1}(0)$ of
$T^*\met(\Sigma)$.  In particular, this means that $\Hcal\subset
(\wp^{-1}(0)\cap p^{-1}(0))/\diff_0(\Sigma)$.

Granted these last points, note next that the quotient of
$\wp^{-1}(0)\cap p^{-1}(0)$ by the semi-direct product of
$\diff_0(\Sigma)$ and $C^{\infty}(\Sigma)$ is $T^*\Tcal$.  In
addition, the pull-back of the form in~\eqref{4.1} to $\wp^{-1}(0)\cap
p^{-1}(0)$ is the pull-back from $T^*\Tcal$ of the canonical cotangent
bundle symplectic form.  

Coupled with the remarks of the previous paragraph, this then has the
following consequence: The form in~\eqref{4.2} is non-degenerate on
$\Hcal$ at non-critical points of the projection to $T^*\Tcal$.  Even
so, it must annihilate the tangent vectors to $\Hcal$ that are mapped
to $0$ by the differential of the projection to $T^*\Tcal$.

Meanwhile, the smooth locus in $\M$ also has a canonical symplectic form;
this induced from the complex, symplectic pairing on
$C^{\infty}(T^*\Sigma\otimes E)$ that
assigns to an ordered pair $(v, v')$ of sections the number
\begin{equation}\label{4.4}
  \int_{\Sigma}\Tr(v\wedge v').
\end{equation}
Here, $\Tr$ denotes the contraction of indices using the $\C$--bilinear inner
product on $E$ that is induced by the product metric on $T^*\Sigma\oplus\R$.  

To elaborate, note that when $\nabla$ is the covariant derivative of a
flat, irreducible $\so{3}{\C}$ connection on $E$, then the tangent space at
the image of $\nabla$ in $\M$ can be identified with the quotient of the kernel
of the corresponding covariant exterior derivative,
\begin{equation}
  d^{\nabla}\co C^{\infty}(T^*\Sigma\otimes E) \to%
  C^{\infty}(\Lambda^2T^*\Sigma\otimes E)
\end{equation}
by the image of $\nabla$ from $C^{\infty}(E)$.  This granted, an
appeal to Stokes' theorem explains why the symplectic pairing
in~\eqref{4.4} sends a pair $(v, v')$ to zero whenever
$v\in\operatorname{kernel}(d^{\nabla})$ and
$v'\in\operatorname{image}(\nabla)$.  Thus, \eqref{4.4} descends as
a closed form to $\M$.  The desired symplectic form is obtained by
first taking the imaginary part of the resulting complex-valued form
and then multiplying the result by the seemingly perverse factor of
$-2\sqrt{6}$.  Denote this form by $\omega_{\M}$.  To see that
$\omega_{\M}$ is non-degenerate, recall first that the quotient of the
kernel of $d^{\nabla}$ by the image of $\nabla$ is isomorphic to the
vector space of sections $v$ of $T^*\Sigma\otimes E$ that obey both
$\epsilon_{AB}\nabla_Av_B=0$ and also $\overline{\nabla}_Av_A=0$.
This vector space admits the anti-holomorphic involution that sends $v$
to $t(v)$ where $t(v)_A\equiv i\epsilon_{AB}\overline{v_B}$.  Then
$\omega_{\M}(t(v),v) = \int_{\Sigma}|v|^2d\operatorname{vol}$.

Now, the map from $\Hcal$ to $\M$ pulls back the form $\omega_{\M}$ to
a closed form on $\Hcal$ and one has:

\begin{lem}\label{lem:4.1}
The pull-backs to $\Hcal$ of the symplectic forms on $T^*\Tcal$ and $\M$ agree.
\end{lem}

\begin{proof}[Proof of Lemma~\ref{lem:4.1}]
Let $(g, m)\in\Hcal$ and let $\nabla$ denote the resulting flat
$\so{3}{\C}$ connection from~\eqref{1.8} and~\eqref{1.9}.  When $(h,
n)$ is a tangent vector to $\Hcal$ at $(g, m)$, then the differential
of the map from $(g, m)$ to $\M$ assigns the $d^{\nabla}$--closed
vector $v$ in~\eqref{3.4} to $(h, n)$.  If $(h, n)$ and $(h',
n')$ are two such vectors, then the imaginary part of the pairing
in~\eqref{4.4} on the corresponding ordered pair $(v, v')$ is
exactly $(-2\sqrt{6})^{-1}$ times the pairing in~\eqref{4.3} between
$(h,n)$ and $(h', n')$.
\end{proof}

This last result has the following amusing corollary:

\begin{lem}\label{lem:4.2}
Let $(g, m)\in\Hcal$ and let $K_{\Tcal}$ and $K_{\M}$ denote the
respective kernels in $T\Hcal|_{(g,m)}$ of the differentials of the maps to
$T^*\Tcal$ and to $\M$.  Then the following are true:
\begin{itemize}
\item $K_{\Tcal}\cap K_{\M} = \{0\}$ and so the differential of the
      map to $T^*\Tcal$ is injective on $K_{\M}$ while that of the map to $\M$
      is injective on $K_{\Tcal}$.
\item The cokernel at the image of $(g, m)$ in $T^*\Tcal$ of the
      differential of the map from $\Hcal$ is the symplectic dual to
      the image of $K_{\M}$, while the cokernel at the image of $(g,
      m)$ in $\M$ of the map from $\Hcal$ is the symplectic dual of the
      image of $K_{\Tcal}$.
\item With suitable signs chosen for the symplectic forms on
      $T^*\Tcal$ and $\M$, the maps to these spaces immerse $\Hcal$ in
      $T^*\Tcal\times\M$ as a immersed, Lagrangian subvariety.
\end{itemize}
\end{lem}

With regards to this lemma, keep in mind that
Propositions~\ref{prop:3.2} and~\ref{prop:3.3} respectively identify
$K_{\M}$ and $K_{\Tcal}$ with the kernel of the operator $\Delta\equiv-d_Cd_C+(\frac{1}{3}-|m|^2)$.

\begin{proof}[Proof of Lemma~\ref{lem:4.2}]
To prove that $K_{\Tcal}\cap K_{\M} = \{0\}$, appeal to
Lemma~\ref{lem:D.1} in the Appendix to justify the remark that any
element in this intersection provides a non-zero pair of sections
$(\eta, u)$ of $E$ that obey $\nabla_A\eta = \epsilon_{AB}\nabla_B'u$
As argued in Section~\ref{section:structure}, there are no such pairs.

Granted the first point of Lemma~\ref{lem:4.1}, here is why the second
point holds: Since $K_{\M}$ is annihilated by the differential of the
map to $\M$, it must have zero pairing with the whole of
$T\Hcal|_{(g,m)}$ using the pull-back of the symplectic from from $\M$.
Thus, it has zero pairing with respect to the pull-back of the form
from $T^*\Tcal$.  However, as it is mapped injectively to $T(T^*\Tcal)$
by the differential of the map from $\Hcal$, this can happen only if the
symplectic dual of its image in $T(T^*\Tcal)$ has trivial intersection
with the image of $T\Hcal|_{(g,m)}$.  A comparison of dimensions then
establishes the assertion that the symplectic dual to the image of $K_{\M}$
in $T(T^*\Tcal)$ is isomorphic to the cokernel of the differential at
$(g, m)$ of the map to $T^*\Tcal$.  The analogous argument where the roles
of $\M$ and $T^*\Tcal$ are switched proves the assertion in the second
point of the lemma about $K_{\Tcal}$ and the cokernel of the differential of
the map to $\M$.

The third point of the lemma follows as an immediate consequence of
Lemma~\ref{lem:4.1} and the second point.  In this regard, note that
this consequence was pointed out to the author by Curt McMullen.
\end{proof}

The nullity of an immersed, minimal surface in a three manifold is
defined to be the dimension of the null space of the Hessian at zero
of a certain function on the space of sections of the normal bundle.
To be precise, the function assigns to any given section the area of a
corresponding surface in the $3$--manifold; and the corresponding
surface is obtained by moving the original surface distance $1$ along
the geodesics that are tangent to the given section.  When the ambient
$3$--manifold is hyperbolic, the null space of this Hessian is
precisely the kernel of $-d_Cd_C + (\frac{1}{3} - |m|^2)$.  This
understood, then Lemma~\ref{lem:4.2} leads to

\begin{lem}\label{lem:4.3}
Let $M$ be any hyperbolic $3$--manifold and let $\Sigma\subset M$ be a compact,
oriented, immersed minimal surface with Euler characteristic $-\chi$.  Then
the nullity of $\Sigma$ is no greater than $3\chi$.
\end{lem}

\section{Extending the germ as an honest hyperbolic metric}%
\label{section:extending}
The purpose of this section is to construct from any given $(g, m)\in
\Hcal$ a hyperbolic metric (with scalar curvature $-1$) on a
neighborhood of $\Sigma\times\{0\}$ in $\Sigma\times\R$ whose
respective first and second fundamental forms on $\Sigma\times\{0\}$
are $g$ and $m$.  In this regard, the model for such an extension is
that of the pair $(g, 0)$ where $g$ is a metric on $\Sigma$ with
constant scalar curvature equal to $-\frac{1}{3}$.  In particular, the
latter pair comes from a complete hyperbolic metric on
$\Sigma\times\R$, this the metric given on the product of $\R$ with a
local coordinate chart on $\Sigma$ by the line element with square
\begin{equation}
  ds^2=\cosh^2(\frac{1}{\sqrt{6}}t)\cdot g_{AB}(z)dz^Adz^B+dt^2.
\end{equation}
Now, suppose that $(g, m)$ is some given pair from $\Hcal$.  The plan
is to seek a hyperbolic metric near $\Sigma\times \{0\}$ whose line element has
square
\begin{equation}\label{5.2}
  ds^2=\gamma_{AB}(t,z)dz^Adz^B+dt^2.
\end{equation}
As it turns out, such a metric can be obtained by solving at each
point in $\Sigma$ the ordinary differential equation
\begin{equation}\label{5.3}
  \begin{split}
    \partial_t\gamma & = 2\gamma\mu\\
    \partial_t\mu & = -\mu^2+\frac{1}{6}\mathbb{I}.
  \end{split}
\end{equation}
where $\gamma$ and $\mu$ are $2\times 2$ matrix functions of $t$ with
$t = 0$ values $\gamma = g$ and $\mu = g^{-1}m$.  In this regard,
standard techniques from the theory of ordinary differential equations
can be used to prove that~\eqref{5.3} with the attending $t = 0$ conditions
has a unique solution on some interval about $t = 0$.  Moreover, the
size of this interval can be bounded from below in terms of $g$ and $m$.
This last point implies that there is a fixed interval about $0$ in $\R$ on
which~\eqref{5.3} and its attending $t = 0$ conditions has a unique solution
for each point in $\Sigma$.  

The equations of Gauss--Codazzi can be used to verify that the
resulting $\gamma(t, z)$ defines a metric via~\eqref{5.2} on a
neighborhood of $\Sigma\times\{0\}$ in $\Sigma\times\R$ that is
hyperbolic with scalar curvature $-1$.

\section{A neighborhood of the Fuchsian locus}\label{section:neighborhood}
The purpose of this section is to describe the structure of $\Hcal$
near the locus of pairs of the form $(g_0, 0)$ where $g_0$ is a metric
on $\Sigma$ with constant curvature equal to $-\frac{1}{3}$.  This is
called the ``Fuchsian locus'' because the image such a point in $\M$
is a Fuchsian group.  In any event, this locus is observedly a smooth
submanifold of $\Hcal$.  Moreover, the latter has a tubular
neighborhood that consists of the pairs $(g, m)$ that obey (1.1) with
the auxiliary bound
\begin{equation}\label{6.1}
  |m|^2 < \frac{1}{3}.
\end{equation}
In fact, the bound in~\eqref{6.1} insures that the operator
$\Delta_{(g,m)} = -d_Cd_C - |m|^2 + \frac{1}{3}$ is strictly positive,
and so the map from $\Hcal$ to $T^*\Tcal$ identifies the open subset
of $\Hcal$ where~\eqref{6.1} holds with an open neighborhood of the
zero section in $T^*\Tcal$.  Let $\U$ denote the subset in $\Hcal$
where~\eqref{6.1} holds.

The set $\U$ has the following property:

\begin{prop}\label{prop:6.1}
Each $(g, m)\in\U$ consists of the respective first and second
fundamental form of a surface in $\R\times\Sigma$ that is minimal with
respect to some complete hyperbolic metric.  Meanwhile, the
corresponding flat $\so{3}{\C}$ connection defines a homomorphism from
$\pi_1(\Sigma)$ to $\psl{2}{\C}$ whose image is a quasi-Fuchsian group.
\end{prop}

\begin{proof}[Proof of Proposition~\ref{prop:6.1}]
As explained in the previous section, the equations in~\eqref{5.3}
have a unique solution on an open neighborhood in $\Sigma\times\R$ of
$\Sigma\times\{0\}$ if $\gamma|_{t=0} = g$ and $\mu|_{t=0} = g^{-1}m$.
This noted, suppose that this solution to~\eqref{5.3} extends to the
whole of $\Sigma\times\R$.  Granted that such is the case, it then
follows from the algebraic structure of~\eqref{5.3} that $\gamma$ is a
positive definite, symmetric section over $\Sigma\times\R$ of
$T^*\Sigma\otimes T^*\Sigma$ and thus~\eqref{5.2} defines a complete,
hyperbolic metric on the whole of $\Sigma\times\R$.  Of course, by
virtue of the $t = 0$ conditions, the surface $\Sigma\times\{0\}$ is
minimal with respect to this metric with $(g, m)$ being its respective
first and second fundamental forms.

To prove that the solution does indeed extend to $\Sigma\times\R$, let $x$ denote the trace of $\mu$ and let $y$ denote the trace of $\mu^2$.  Since 
\begin{equation}\label{6.2}
  \tr(\mu^3) =\frac{1}{2}(3\tr(\mu)\tr(\mu^2)-(\tr\mu)^3),
\end{equation}
the second equation implies that
\begin{equation}\label{6.3}
  \begin{split}
    \partial_tx &= -y+\frac{1}{3},\quad\text{and}\\
    \partial_ty &= -x(3y-x^2-\frac{1}{3}).
  \end{split}
\end{equation}
Of interest are the solutions to these equations that start at $t = 0$
with $x = 0$ and $y < \frac{1}{3}$.  Take note that the given initial
conditions also imply that $\frac{1}{3}x^2\leq y$ at $t = 0$.  

The analysis that follows considers the $t\geq 0$ evolution.  The
identical analysis holds for $t\leq 0$ after reversing the sign of $\mu$.
In any event, to study the $t\geq 0$ case, note first that the compact
subset where
\begin{equation}\label{6.4}
  \begin{split}
    0 &\leq x\quad\text{and}\\
    \frac{1}{2}x^2&\leq y\leq\frac{1}{3}
  \end{split}
\end{equation}
is a trapping region for the trajectories that obey~\eqref{6.3}.  This
is to say that no trajectory can exit this region.  To verify this
claim, simply check the direction of motion implied by~\eqref{6.3} on
the boundary.  In this regard, note that
\begin{equation}\label{6.5}
  \partial_t(y-\frac{1}{2}x^2) = -2x(y-\frac{1}{2}x^2).
\end{equation}
The fact that the no trajectory can leave the region where~\eqref{6.4} holds
implies that all trajectories that start in this region are defined
for all $t\geq 0$.  

The next point to make is that the flow in~\eqref{6.3} has two fixed
points, both at corners of the region that is delineated
by~\eqref{6.4}; one fixed point is $(x = 0, y = \frac{1}{3})$ and the
other is the point $(x = \sqrt{\frac{2}{3}}, y = \frac{1}{3})$. Since
$\partial_tx > 0$ where $y <\frac{1}{3}$ , the first point is
repelling where $y <\frac{1}{3}$ .  On the other hand, the second
point is an attracting fixed point.  Moreover, as $\partial_tx > 0$
where $y <\frac{1}{3}$ , all trajectories that start
where~\eqref{6.4} holds and with $y <\frac{1}{3}$ limit as
$t\rightarrow\infty$ to the $(x = \sqrt{\frac{2}{3}}, y =
\frac{1}{3})$ fixed point.

To see what this all implies for solutions to~\eqref{6.3}, introduce
$\nu\equiv \mu-\frac{1}{2}x\mathbb{I}$.  Thus, $\nu$ is the traceless
part of $\mu$.  Then~\eqref{5.3} asserts that
\begin{equation}\label{6.6}
  \partial_t\nu = -x\nu.
\end{equation}
Since $x\in (0, \sqrt{\frac{2}{3}})$ for all $t\geq 0$, this last
equation implies that the solution $\nu(t)$ to~\eqref{6.6} exists for
all $t\geq 0$. Moreover,
\begin{equation}\label{6.7}
  |\nu(t)|\leq \text{constant}\,e^{-\sqrt{\frac{2}{3}}t}.
\end{equation}
Thus, the solution $\mu$ to the second line in~\eqref{5.3} exists for
all $t\geq 0$.

Granted~\eqref{6.7}, it then follows from the first line
in~\eqref{5.3} that $\gamma$ exists for all $t > 0$ as well.  Moreover,
\begin{equation}
  \gamma\sim\gamma_+e^{\sqrt{\frac{2}{3}}t}
\end{equation}
as $t\rightarrow\infty$, where $\gamma_+$ is a smooth metric on
$\Sigma$.  

With the preceding understood, it remains only to establish that the
image of the pair $(g, m)$ in $\M$ defines a quasi-Fuchsian
representation of $\pi_1(\Sigma)$ into $\psl{2}{\C}$.  To this end,
note that $\Sigma\times\R$ with the $(g, m)$ version of the metric
in~\eqref{5.2} is uniformized by the hyperbolic ball, $\mathbb{B}$,
with its metric of constant sectional curvature $-\frac{1}{3}$.  This
is to say that the space $\Sigma\times\R$ with its hyperbolic metric
is the quotient of $\mathbb{B}$ by the image, $\Gamma$, of
$\pi_1(\Sigma)$ in $\psl{2}{\C}$.  The inverse image of
$\Sigma\times\{0\}$ via this map is a properly embedded disk in
$\mathbb{B}$ with principle curvatures having absolute value less than
$\frac{1}{\sqrt{6}}$.  This understood, theorems of Epstein
in~\cite{epstein1} and~\cite{epstein2} imply that $\Gamma$ is quasi-Fuchsian.
\end{proof}

\section*{Appendices}\addcontentsline{toc}{secion}{Appendices}
\appendix\begingroup\small
\section{The curvature for the metric in~\eqref{eq:1.2}}
The computation begins with the definition of the orthonormal frame 
\begin{equation}
  \begin{split}
    e^A & = dz^A+tm_{AB}dz^B+\frac{1}{12}t^2dz^A,\quad\text{and}\\
    e^3 & = dt,
  \end{split}
\end{equation}
where $(z^A)$ are Gaussian normal coordinates on $\Sigma$ for the metric $g$. A
computation then finds
\begin{equation}
  \begin{split}
    de^A & = td_Cm_{AB}dz^C\wedge dz^B+m_{AB}dt\wedge
    dz^B+\frac{1}{6}tdt\wedge dz^A\\
    & = -\Gamma^{AB}e^B-\Gamma^{A3}e^3\\
    de^3 &= 0 = -\Gamma^{A3}\wedge e^A.
  \end{split}
\end{equation}
Thus, 
\begin{equation}\label{A.3}
  \begin{split}
    \Gamma^{AB} & = t(d_Bm_{AC}-d_Am_{BC})dz^C,\quad\text{and}\\
    \Gamma^{A3} & = m_{AB}dz^B + \frac{1}{6}tdz^A.
  \end{split}
\end{equation}
To continue, the computation of the curvature $2$--form
$\mathcal{R}^{ij} = \frac{1}{2}R^{ij}_{km} e^k\wedge e^m$ finds it at
$t = 0$ equal to
\begin{equation}
  \begin{split}
    \mathcal{R}^{AB} &=
    \frac{1}{4}r\epsilon^{AB}\epsilon_{CD}dz^Cdz^D
    +(d_Bm_{AC}-d_Am_{BC})dt\wedge dz^C\\
    &\quad\quad\quad - m_{AC}m_{BD}dz^C\wedge dz^D,\\
    \mathcal{R}^{A3} &= d_Cm_{AB}dz^C\wedge dz^B+\frac{1}{6}dt\wedge dz^A.
  \end{split}
\end{equation}
Thus, the curvature form at $t = 0$ is
\begin{equation}
  \begin{split}
    R^{AB}_{CD} & =
    \frac{1}{2}r\epsilon^{AB}\epsilon_{CD}-(m_{AC}m_{BD}-m_{AD}m_{BC}),\\
    R^{AB}_{C3} & = -(d_Bm_{AC}-d_Am_{BC}),\\
    R^{A3}_{B3} & = -\frac{1}{6}g_{AB}.
  \end{split}
\end{equation}
This gives the $t = 0$ Ricci tensor with components
\begin{equation}
  \begin{split}
    R_{AC} & =
    \frac{1}{2}(r-\frac{1}{3})g_{AC}+m^0_{AC}m^0_{BC}-\frac{1}{4}k^2g_{AB},\\
    R_{A3} & = d_Bm^0_{AB}-\frac{1}{2}d_Ak,\\
    R_{33} & = -\frac{1}{3},
  \end{split}
\end{equation}
where $m^0$ is the traceless part of $m$.  Thus, the $t = 0$ Ricci tensor
obeys the hyperbolicity condition if and only if
\begin{equation}
    r = -(\frac{1}{3}+|m^0|^2)+\frac{1}{2}k^2\quad\text{and}\quad%
    d_Cm_{AB}-d_Bm_{AC} = 0.
\end{equation}
Note that the second fundamental form for the $t = 0$ slice can be
computed as follows:
\begin{equation}
  \begin{split}
    \kappa_{AB} & = \iproduct{\nabla_Ae^3}{dz_B}\\
     & = -\Gamma^{3B}_A\\
     & = m_{AB},
  \end{split}
\end{equation}
by virtue of~\eqref{A.3}.  Thus, the $t = 0$ slice has zero mean
curvature when $k = g^{AB}m_{AB} = 0$.

\section{The variation of (1.1) with a conformal change of the metric}
The differential of the equation of (1.1) at a pair $(g, m)$ can be
computed in the following manner: First of all, consider the change
where $(g, m)\rightarrow (e^{-u}g, m)$.  The change in the covariant
derivative is computed as follows: Observe first that the new
orthonormal frame is $\{\hat{e}^A = e^{-\frac{u}{2}}e^A\}$ and
\begin{equation}
  d\hat{e}^A = -\frac{1}{2}d_Cue^{-\frac{u}{2}}e^C\wedge
  e^A-\gamma^{AB}\wedge \hat{e}^B.
\end{equation}
Thus, 
\begin{equation}
  \hat{\gamma}^{AB} = \gamma^{AB} - \frac{1}{2}(e^Ad_Bu-e^Bd_Au).
\end{equation}
This then finds the new curvature $2$--form
\begin{equation}
  \hat{r}^{AB} = r^{AB} - \frac{1}{2}(d_Cd_Bue^C\wedge e^A -
  d_Cd_Aue^C\wedge e^B).
\end{equation}
Hence,
\begin{equation}
  \hat{r}^{AB}_{CD} = e^u(r^{AB}_{CD} - \frac{1}{2}(d_Cd_Bug_{AD}-d_Dd_Bug_{AC}-d_Cd_Aug_{BD}+d_Dd_Aug_{BC})).
\end{equation}
Taking the necessary traces finds
\begin{equation}\label{B.5}
  \hat{r} = e^u(r+\Delta u).
\end{equation}
In terms of the coordinate frame $\{dz^A\}$, the new Christoffel symbols are
\begin{equation}
  \begin{split}
  \delta\gamma^D_{CB} & \equiv \hat{\gamma}^D_{CB}-\gamma^D_{CB}\\
   & = \frac{1}{2}(d_Bug_{DC}+d_Cug_{DB}-d_Dug_{CB}).
  \end{split}
\end{equation}
This then implies that $\hat{d}_Cm_{AB}-\hat{d}_Am_{CB}$ is equal to
\begin{equation}
  \begin{split}
    =\; & -\delta\gamma^D_{CB}m_{AD}+\delta\gamma^D_{AB}m_{CB}\\ 
    =\;& \frac{1}{2}(d_Bum_{AC}+d_Cum_{AB}-d_Dum_{AD}g_{CB})\\
    &-\frac{1}{2}(d_Bum_{CA}+d_Aum_{CB}-d_Dum_{CD}g_{AB})\\
    =\; &
    \frac{1}{2}\epsilon_{CA}(d_Du\epsilon_{DF}m_{FB}-d_Dum_{DF}\epsilon_{BF})\\
    =\; & 0
  \end{split}
\end{equation}
by virtue of the fact that the tensor $\epsilon_{DF}m_{FB}$ is
symmetric when $m$ is symmetric and tracefree.  Thus, $(e^{-u}g, m)$
obeys the first equation in~\eqref{eq:1.1a} when $(g, m)$ does.
Needless to say, the equation in~\eqref{eq:1.1b} is also obeyed by
$(e^{-u}g, m)$.  

As for the second equation in~\eqref{eq:1.1a}, introduce
\begin{equation}
  f(g,m)\equiv r_g+(|m|^2+\frac{1}{3}).
\end{equation}
Then by virtue of~\eqref{B.5},
\begin{equation}
  \begin{split}
    f(e^{-u}g,m) & = e^u(r+\Delta u+e^u|m|^2+\frac{1}{3}e^{-u})\\
    & = e^u(\Delta u + (e^u-1)|m|^2+\frac{1}{3}(e^{-u}-1)).
  \end{split}
\end{equation}
Note that if $f(e^{-u}g,m)$ is to vanish for small $u$, then
\begin{equation}
  \Delta u + (|m|^2-\frac{1}{3})u+O(u^2) = 0.
\end{equation}
\section{The connection defined by~\eqref{2.5}}
The curvature of the covariant derivative in~\eqref{2.5} at a given
point in $\Sigma$ is computed most easily by choosing the frame
$\{e^A\}$ for $T^*\Sigma$ to have zero metric covariant derivative at
the point in question.  This understood, then the curvature $2$--form
is the $\so{3}{\C}$--valued $2$--form whose components are
$\mathcal{R}_{AC} = [\nabla_A, \nabla_C]$.  Therefore, two appeals
to~\eqref{2.5} find
\begin{equation}
  \mathcal{R}_{AC}\eta =%
   \begin{pmatrix}
    -r_{DBAC}\eta_D+(d_A\theta_{CB}-d_C\theta_{AB})\eta_3 -
    (\theta_{AB}\theta_{CD}-\theta_{CB}\theta_{AD})\eta_D\\
    -(d_A\theta_{CD}-d_C\theta_{AD})\eta_D-(\theta_{AD}\theta_{CD}-\theta_{CD}\theta_{AD})\eta_3
 \end{pmatrix},
\end{equation}
where $r_{DBAC}=\frac{1}{2}\epsilon_{DB}\epsilon_{AC}r$ is the
curvature tensor for the Riemannian metric $g$.  Now, the bottom term
here is zero provided that
\begin{equation}
  d_A\theta_{CD}-d_C\theta_{AD} = 0.
\end{equation}
This noted, the top term is zero provided that
\begin{equation}
  -r-\theta_{AB}\theta_{BA}+\theta_{AA}\theta_{BB} = 0.
\end{equation}
These are the equations in (1.1).

The next task for this appendix is to verify that any two solutions to
(1.1) from the same $\diff_0$ orbit define gauge equivalent flat
$\so{3}{\C}$ connections.  For this purpose, it proves more convenient
to work with the metric in~\eqref{eq:1.2} on $\Sigma\times\R$.  To
start, let $g$ now denote a given metric on $\Sigma\times\R$.  Fix an
oriented, orthonormal frame, $\{e^i\}_{i=1,2,3}$, for the cotangent
bundle of $\Sigma\times\R$ and this frame can then be used to define
from $g$ the $\so{3}{\C}$ connection on the trivial $\C^3$ whose covariant
derivative, $\nabla$, is
\begin{equation}\label{C.4}
  (\nabla\eta)_i = d\eta_i-\Gamma^k_i\eta_k+\frac{i}{\sqrt{6}}\epsilon_{ijk}\eta_ke^j.
\end{equation}
Here, $\Gamma^k_i$ denotes the Levi--Civita connection as defined by
the given metric using the given orthonormal frame.  The covariant
derivative defined by~\eqref{2.5} in directions along the $t = 0$
slice is the restriction of the version of~\eqref{C.4} that is defined
by the metric in~\eqref{eq:1.2}.  If a first order change in the
metric $g$ is given by $\delta g = h$ with $h$ a symmetric, $3\times 3$
tensor, then the associated orthonormal frame changes to first order
so that
\begin{equation}
  \delta e^k = \frac{1}{2} h_{kj} e^j.
\end{equation}
Meanwhile, the metric's Levi--Civita covariant derivative changes to
first order as
\begin{equation}
  \begin{split}
  \delta\Gamma^k_i &=
  \frac{1}{2}(\partial_ih_{kj}-\partial_kh_{ij})e^j\\
   & = \epsilon_{imk}(\frac{1}{2}\epsilon_{nmp}\partial_nh_{pj}e^j)
  \end{split}
\end{equation}
Here, $\partial$ is used to denote the covariant derivative of the
metric $g$.  Thus, the covariant derivative in~\eqref{C.4} changes to first
order as
\begin{equation}
  \delta(\nabla\eta)_i = \epsilon_{imk}(-\frac{1}{2}\epsilon_{nmp}\partial_nh_{pj}+\frac{i}{\sqrt{6}}h_{mj})e^j\eta_k.
\end{equation}
As a parenthetical remark, note that such a change preserves the
flatness condition to first order if and only if the $\C^3$--valued $1$--form
\begin{equation}\label{C.9}
  v_i\equiv (-\frac{1}{2}\epsilon_{nip}\partial_nh_{pj}+\frac{i}{\sqrt{6}}h_{ij})e^j
\end{equation}
obeys $D^{\nabla}v = 0$ where $D^{\nabla}$ denotes the exterior
covariant derivative that is defined by $\nabla$.  

By comparison, a change in the covariant derivative comes from a first
order change in the connection in a direction tangent to its gauge
orbit if and only if the change in $\nabla$ has the form
\begin{equation}\label{C.10}
  \delta(\nabla\eta)_i = \epsilon_{imk}(\nabla\sigma)_m\eta_k
\end{equation}
where $\sigma$ is some $\C^3$--valued function.

Granted the preceding, suppose now that the variation, $h$, in the
metric $g$ comes from the action by Lie derivative of a vector field,
$v$, on $\Sigma\times\R$.  This being the case, then $h_{ij} =
\partial_iv_j + \partial_jv_i$ and so~\eqref{C.9} reads
\begin{equation}
  \begin{split}
    v_i & =
    -\frac{1}{2}\epsilon_{nip}\partial_n(\partial_pv_j+\partial_jv_p)e^j
    + \frac{i}{\sqrt{6}}(\partial_iv_j+\partial_jv_i)e^j\\
    & = \nabla(-\frac{1}{2}\epsilon_{nip}\partial_nv_p+\frac{i}{\sqrt{6}}v_i).
  \end{split}
\end{equation}
As this has the form given in~\eqref{C.10}, the infinitesimal action
of $\diff(\Sigma\times\R)$ on the metric results in an infinitesimal
action of the gauge group on the resulting connection.

To finish the story, note that this first order result can then by
``integrated'' to establish that any two solutions to (1.1) that lie
on the same $\diff_0(\Sigma)$ orbit define gauge equivalent, flat
$\so{3}{\C}$--connections along $\Sigma$.
 
\section{The cokernel of the operator $L_{(g,m)}$}\label{appendix:D}
A pair $(\sigma_B, \sigma_3)$ is in the cokernel of $L_{(g,m)}$ if and
only if it is $L^2$--orthogonal to all pairs $(\gamma_B, \gamma_3)$
that can be written as in~\eqref{2.1}.  This then implies that the pair
satisfy the following system of equations:
\begin{gather}\label{D.1}
    -\frac{1}{2}(\epsilon^{CA}d_C\sigma_B+\epsilon^{CB}d_C\sigma_A-g_{AB}\epsilon_{EF}d_E\sigma_F)+2m_{AB}\sigma_3
    = 0\\
   \begin{align*}
     & \left(\frac{1}{2}(\frac{1}{3}-|m|^2)\sigma_3-d_Cd_C\sigma_3\right)g_{AB}
     +d_Ad_B\sigma_3-\frac{1}{4}\epsilon^{EC}d_C(\sigma_Bm_{EA}+\sigma_Am_{EB})\\ 
    &\quad-\frac{1}{4}\left(\epsilon^{EA}d_E(\sigma_Cm_{CB})+\epsilon^{EB}d_E(\sigma_Cm_{CA})\right)-\frac{1}{2}\epsilon^{CD}d_C\sigma_Dm_{AB}
    = 0.
    \end{align*}   
\end{gather}
Now write $\sigma_B=\epsilon_{BC}v_C$ and $\sigma_3=-\frac{1}{2}v_3$.
As is explained next, the equations just written are equivalent to the
following:
\begin{align}\label{D.2}
  & d_Av_B+d_Bv_A-g_{AB}d_Cv_C+2m_{AB}v_3 = 0\\
  & \frac{1}{2}(\frac{1}{3}+|m|^2)v_3g_{AB}-d_Ad_Bv_3+v^Cd_Cm_{AB}+d_Bv_Cm_{CA}+d_Av_Cm_{CB}
  = 0\notag
\end{align}

To start the explanation, make the indicated substitution and then up
to a constant factor, the following is the result:
\begin{gather}\label{D.3}
  d_Av_B+d_Bv_A-g_{AB}d_Cv_C+2m_{AB}v_3 = 0\\
  \begin{align*}
   & \left(d_Cd_Cv_3-\frac{1}{2}(\frac{1}{3}-|m|^2)v_3\right)g_{AB}-d_Ad_Bv_3\\
   &\quad-\frac{1}{2}(2d_Cv_Cm_{AB}-d_Bv_Cm_{CA}-d_Av_Cm_{CB})\\
   &\quad-\frac{1}{2}(d_Cv_Am_{CB}+d_Cv_Bm_{CA}-2d_Cv_Cm_{AB}-2v^Cd_Cm_{AB})\\
   &\quad+d_Dv_Dm_{AB} = 0.
  \end{align*}
\end{gather}
Note that the first equation in~\eqref{D.2} is identical to its mate
in~\eqref{D.3}.  To simplify the second equation, begin by taking its
trace to deduce that
\begin{equation}\label{D.4}
  d_Cd_Cv_3-(\frac{1}{3}-|m|^2)v_3 = 0.
\end{equation}
The second equation in~\eqref{D.3} can now be rewritten
using~\eqref{D.4} and the first equation in~\eqref{D.3} to read
\begin{multline}\label{D.5}
  \frac{1}{2}(\frac{1}{3}-|m|^2)v_3g_{AB}-d_Ad_Bv_3-\frac{1}{2}(2d_Cv_Cm_{AB}-2d_Bv_Cm_{CA}-2d_Av_Cm_{CB})\\
  -\frac{1}{2}\Big((g_{CA}d_Dv_D-2m_{CA}v_3)m_{CB}+(g_{CB}d_Dv_D-2m_{CB}v_3)m_{CA}\\
  -2d_Cv_Cm_{AB}-2v^Cd_Cm_{AB}\Big) +d_Dv_Dm_{AB} = 0.
\end{multline}
The elimination of canceling terms and the use of the identity
$m_{AB}m_{BC}=\frac{1}{2}|m|^2g_{AC}$ makes~\eqref{D.5} into the
second equation in~\eqref{D.2}.

The next part of this appendix explains why~\eqref{D.2} is equivalent
to~\eqref{2.6}.  As the next lemma points out, a more general version
of~\eqref{D.2} and~\eqref{2.5} together imply~\eqref{2.6}.

\begin{lem}\label{lem:D.1}
Suppose that $(v_A, v_3)$ define a real-valued section of
$T^*\Sigma\oplus\R$, that $\kappa$  is a real-valued constant, and that the
following equations hold:
\begin{equation}\label{D.6}
  \begin{split}
  & d_Av_B+d_Bv_A-g_{AB}(d_Cv_C+\kappa)+2m_{AB}v_3 = 0\\
  & \frac{1}{2}(\frac{1}{3}+|m|^2)v_3g_{AB}-d_Ad_Bv_3+v^Cd_Cm_{AB}+d_Bv_Cm_{CA}+d_Av_Cm_{CB}
  = 0.
  \end{split}
\end{equation}
Then the complex-valued sections $\eta$ and $u$ of
$(T^*\Sigma\oplus\R)_{\C}$ given by
\begin{equation}\label{D.7}
  \begin{split}
    & \eta_B\equiv
    v_B+i\sqrt{6}(-\epsilon_{CB}d_Cv_3+v_E\epsilon_{CB}m_{EC}),\\
    & \eta_3\equiv v_3+i\sqrt{\frac{3}{2}}\epsilon_{EF}d_Ev_F,\\
    & u_B\equiv 0,\quad\text{and}\\
    & v_3 = -i\sqrt{\frac{3}{2}}(d_Cv_C+\kappa)
  \end{split}
\end{equation}
obey the equation $\nabla_A\eta=\epsilon_{AB}\nabla_B'u$
\end{lem}

\begin{proof}[Proof of Lemma~\ref{lem:D.1}]
Note first that the vanishing of the first term in~\eqref{D.6} implies that
\begin{equation}\label{D.8}
  d_Av_B-\frac{1}{2}\epsilon_{AB}\epsilon_{EF}d_Ev_F+m_{AB}v_3 = \frac{1}{2}g_{AB}(d_Cv_C+\kappa).
\end{equation}
Then, since 
\begin{equation}
  \frac{1}{2}g_{AB}(d_Cv_C+\kappa) = \epsilon_{AC}(-\frac{i}{\sqrt{6}}\epsilon_{CB}u_3),
\end{equation}
this last equation is the real part of the equation
$d_A\eta_B+\theta_{AB}\eta_3 = \epsilon_{AC}(-\theta_{CB}u_3)$.  

To continue, the vanishing of the second term
in~\eqref{D.6} can be rewritten using~\eqref{D.8} so as to read
\begin{multline}
  -d_Ad_Bv_3+\frac{1}{2}(\frac{1}{3}-|m|^2)g_{AB}v_3+v_Cd_Cm_{AB}\\
  +m_{AC}\epsilon_{BC}(\epsilon_{EF}d_Ev_F)
  = -m_{AB}(d_Cv_C+\kappa).
\end{multline}
Moreover, this is equivalent to
\begin{equation}\label{D.11}
  -d_Ad_Bv_3+\frac{1}{6}g_{AB}v_3+d_A(v_Cm_{CB})+\frac{1}{2}m_{AC}\epsilon_{BC}(\epsilon_{EF}d_Ev_F)
  = -\frac{1}{2}m_{AB}(d_Cv_C+\kappa).
\end{equation}
Finally, the latter gives
\begin{multline}
  d_A\sqrt{6}(-\epsilon_{BD}d_Bv_3+v_C\epsilon_{BD}m_{CB})+\frac{1}{\sqrt{6}}\epsilon_{AD}v_3+m_{AD}(\sqrt{\frac{3}{2}}\epsilon_{EF}d_Ev_F)\\
  = -\sqrt{\frac{3}{2}}\epsilon_{BD}m_{AB}(d_Cv_C+\kappa).
\end{multline}
The latter equation is the imaginary part of the equation
$$d_A\eta_D+\theta_{AD}\eta_3=\epsilon_{AC}(-\theta_{CD}u_3).$$
As will now be explained, it is also the case that
$d_A\eta_3-\theta_{AB}\eta_B=\epsilon_{AC}d_Cu_3$. To this end,
consider that the real part of $\theta_{AB}\eta_B$ can be written as
follows:
\begin{align}\label{D.13}
  \re(\theta_{AB}\eta_B) & =
  m_{AB}v_B+\epsilon_{AB}\epsilon_{CB}d_Cv_3-\epsilon_{AB}\epsilon_{CB}m_{EC}v_E\\
  & = d_Av_3\notag\\
  & = \re(d_A\eta_3).\notag
\end{align}
Meanwhile, the imaginary part of $\theta_{AB}\eta_B$ can be written to read
\begin{align}
  \im(\theta_{AB}\eta_B) & =
  \frac{1}{\sqrt{6}}\epsilon_{AB}v_B+\sqrt{6}(-m_{AB}\epsilon_{CB}d_Cv_3+m_{AB}\epsilon_{CB}v_Em_{EC})\\
  & =
  \epsilon_{AB}\left(\frac{1}{\sqrt{6}}v_B+\sqrt{\frac{3}{2}}|m|^2v_B-\sqrt{6}d_C(m_{CB}v_3)\right)\notag\\
  & = -\sqrt{\frac{3}{2}}\epsilon_{AB}\left(rv_B+d_C(2m_{CB}v_3)\right).\notag
\end{align}
On the other hand, commuting derivatives and invoking (D.8) allows the
imaginary part of $d_A\eta_3$ to be written as
\begin{align}
  \im(d_A\eta_3) & = \sqrt{\frac{3}{2}}\epsilon_{EF}(d_Ad_Ev_F)\\
  & =
  \sqrt{\frac{3}{2}}(\epsilon_{EA}\frac{1}{2}rv_E+\epsilon_{EF}d_Ed_Av_F)\notag\\
  & =
  \sqrt{\frac{3}{2}}\left(\epsilon_{BA}\frac{1}{2}rv_B-\epsilon_{EF}d_Ed_Fv_A-\epsilon_{EF}d_E(2m_{AF}v_3)+\epsilon_{EA}d_Ed_Cv_C\right)\notag\\
  & = -\sqrt{\frac{3}{2}}\epsilon_{AB}\left(rv+d_C(2m_{BC}v_3)) + \epsilon_{AE}d_E(-(d_Cv_C+\kappa)\right).\notag
\end{align}
Together, \eqref{D.11} and~\eqref{D.13} imply the claim that $d_A\eta_3-\theta_{AB}\eta_B=\epsilon_{AC}d_Cu_3$.
\end{proof}\endgroup

\subsection*{Acknowledgements}
\addcontentsline{toc}{subsection}{Acknowledgements} 
The author gratefully thanks Hyam Rubinstein and Curt McMullen for
their suggestions, advice and patience.  This work was supported in
part by a grant from the National Science Foundation.

\Addresses\recd

\end{document}